\documentclass[a4paper,11pt]{article}

\usepackage{graphics,float}
\usepackage[latin1]{inputenc} 
\usepackage[T1]{fontenc} 
\usepackage{lmodern} 
\usepackage[english]{babel} 
\usepackage{tikz}
\usetikzlibrary{arrows,shadows}
\usepackage[affil-it,auth-sc]{authblk}

\usepackage{amsfonts,latexsym,amssymb} 
\usepackage{mathrsfs,amsmath}
\usepackage[all]{xy}
\usepackage{pdfpages}

\usepackage{xcolor} 
\usepackage{color} 

\usepackage[numbers]{natbib} 
\usepackage[squaren, Gray, cdot]{SIunits}

\newlength{\defbaselineskip}
\renewcommand{\hat}{\widehat}

\newcommand{\R}{\ensuremath{\mathbb{R}}}

\renewcommand{\phi}{\varphi}

\newcommand{\vecM}{\mathbf{M}}

\newcommand{\vece}{\mathbf{e}}

\newcommand{\vecn}{\mathbf{n}}

\newcommand{\vecv}{\mathbf{v}}

\newcommand{\vectau}{\boldsymbol\tau}

\newcommand{\ds}{\displaystyle}

\newcommand{\diff}{\sigma^c}

\newenvironment{Remarque} {  \normalfont \begin{flushright} $\rightsquigarrow$ \vline \- \begin{minipage}{10cm} \textbf{brouillon.}}
                          {  \end{minipage} \end{flushright}  }

\newcommand{\comment}[1]{}

    \newtheorem{theorem}{Theorem}[section]

    \newtheorem{proposition}[theorem]{Proposition}
 
    \newtheorem{definition}{Definition}[section]

    \newtheorem{remark}{Remark}[section]
    
\renewcommand{\comment}[1]{ {#1} }
\renewcommand{\comment}[1]{ {#1} }

\renewcommand{\deg}{\mathrm{deg}}

\numberwithin{equation}{section}

\begin{document}
\title{Pattern selection in a biomechanical model for the growth of walled cells}
\renewcommand\Affilfont{\itshape}
\setlength{\affilsep}{1em}
\renewcommand\Authsep{, }
\renewcommand\Authand{ and }
\renewcommand\Authands{ and }
\setcounter{Maxaffil}{2}

\author[1]{\textnormal{Vincent} Calvez}
\author[2]{\textnormal{Laetitia} Giraldi}

\affil[1]{UMPA, \'Ecole normale supérieure de Lyon,  Lyon, France. \smallskip}
\affil[2]{UMA, ENSTA-Paristech, Palaiseau, France. \smallskip}

\maketitle

		\newcount\hh
		\newcount\mm
		\mm=\time
		\hh=\time
		\divide\hh by 60
		\divide\mm by 60
		\multiply\mm by 60
		\mm=-\mm
		\advance\mm by \time
		\def\hhmm{\number\hh:\ifnum\mm<10{}0\fi\number\mm}


\begin{abstract}
In this paper, we analyse a model for the growth of three-dimensional
walled cells. In this model the biomechanical expansion of the cell is 
coupled with the geometry of its wall. We consider that the density of 
building material depends on the curvature of the cell wall, thus 
yielding possible anisotropic growth. The dynamics of the axisymmetric 
cell wall is described by a system of nonlinear PDE including a 
nonlinear convection-diffusion equation coupled with a Poisson equation. 
We develop the linear stability analysis of the spherical symmetric 
configuration in expansion. We identify three critical parameters that 
play a role in the possible instability of the radially symmetric shape, 
namely the degree of nonlinearity of the coupling, the effective 
diffusion of the building material, and the Poisson's ratio of the cell 
wall. We also investigate numerically pattern selection in the nonlinear 
regime. All the results are also obtained for a simpler, but similar,
two-dimensional model.
\end{abstract}




%
%

\section{Introduction}

The physical features that account for the acquisition and maintenance of cellular shape is a current problem in experimental and theoretical biology \cite{chang_how_2014,Dumais13,HaylesNurse01}. 

In this work, we investigate a generic biomechanical model for the growth of walled cells, such as plant cells, fungal hyphae, or fission yeast ({\it S. pombe}). Cell wall can be described as a thin shell subject to a high internal pressure, called the turgor pressure, that can reach up to 10 atmospheres \cite{CampasMahadevan09,Minc09}. Due to this mechanical feature, the cell wall by itself determines the shape of the cell. It is commonly admitted that a good representation of an expanding cell wall is provided by an inflating balloon \cite{Dumais07,chang_how_2014}. 

Several attempts have been made recently to understand the dynamics of growth of walled cells, and in particular rod-like cells and pollen tubes, from a physical and geometrical viewpoint \cite{CampasMahadevan09,DrakeVavylonis13,Dumais06,RojasDumais11}.
Here, we aim to further analyze the coupling between the mechanics of cell wall expansion, and the pattern of growth, the latter being determined by the cell geometry. For this purpose we study a minimal model accounting for the mechanics of wall expansion, and heterogeneous distribution of growth along the wall as a function of its geometry.
%

\bigskip

We follow the works of Dumais et al \cite{Dumais06}, and Drake and Vavylonis \cite{DrakeVavylonis13} for the   biomechanical part of the model. They derive the same set of equations, based on slightly different hypotheses. In the former, the mechanical deformation of the cell wall is described using viscoplasticity theory (see also \cite{CampasMahadevan09} where it is described as a viscous thin shell). The latter develops a model where the wall is viewed as an elastic membrane under pressure, and subject to local remodeling. Apart from the equations,  the common features of these models is the heterogeneity of the mechanical characteristics (either the viscosity in Dumais et al, or the rate of remodeling in Drake and Vavylonis). This heterogeneity (denoted as $\Psi$ in the present article) depends upon the distribution of some growth factor ({\it e.g.} Cdc42 for the fission yeast, wall building polymers for fungal hyphae cells, or cell wall loosening enzymes for plant cells). In this work, we do not attempt to describe with much details the set up of heterogeneity. We shall describe this heterogeneity factor with a generic quasi-stationary reaction-diffusion equation with a geometric source term.

We restrict to the case of an axisymmetric cell. This choice obviously rules out many possible shapes, but it is compatible with the acquisition and maintenance of the rod shape, which is already a challenging problem \cite{chang_how_2014}. Furthermore, the analysis, and even the setting of the mathematical model is dramatically simplified in this context.  
%
The cell boundary is described by the angular deviation of the normal vector from the axis of symmetry, $\varphi(t,s)$, or equivalently by the radius of the cell in the transversal direction, $r(t,s)$, where $s$ is the curvilinear abscissa  (see Figure \ref{Fig:surface_revolution}). The biomechanical model is a system of two equations which determine the velocity of the cell wall in the Frenet frame $\vecv = (v_n, v_{\tau})$, 
\begin{align}
\ds \kappa_{\theta} v_n +\frac{v_{\tau} \cos(\phi)}{r} - \Psi (\sigma_{\theta} -\nu\sigma_s)=0\,, \label{eq:model intro 1} \\ 
\ds \frac{\partial v_{\tau}}{\partial s} - \left( \frac{\kappa_s}{\kappa_{\theta}}\frac{\cos(\phi)}{r}\right) v_{\tau} + \frac{\kappa_s}{\kappa_{\theta}}\Psi\left(\sigma_{\theta}-\nu\sigma_{s}\right)-\Psi \left(\sigma_s-\nu\sigma_{\theta}\right)=0\,. \label{eq:model intro 2}
\end{align}
Equations are derived from the constitutive laws governing the deformation of an elastic cell wall under constant turgor pressure $P$ \cite{DrakeVavylonis13}. The stresses $\sigma_s$ and $\sigma_\theta$ are calculated so as to balance the turgor pressure $P$, 
\[ \sigma_s = \dfrac{P}{2\delta \kappa_\theta}\, , \quad \sigma_\theta = \dfrac{P}{2\delta \kappa_\theta}\left( 2 - \dfrac{\kappa_s}{\kappa_\theta}\right)\, , \]
where $\kappa_s$ and $\kappa_\theta$ denote the principal curvatures, and $\delta$ is the thickness of the wall. As usually done, the latter is assumed to be constant. We aggregate in a single parameter $\Psi$ the (inhomogeneous) mechanical properties of the wall, and the local remodeling rate. 
We will also refer to $\Psi$ as the cell wall extensibility, in reference to Dumais et al \cite{Dumais06}. 
Our main modelling hypothesis is that $\Psi$ depends on the local geometry of the cell, {\it e.g.} via its curvature. More precisely, we assume that $\Psi = F(\mu)$, where $F$ is a certain nonlinear function, and $\mu$ represents the distribution of growth signal \cite{DrakeVavylonis13} at the cell wall. Equation for $\mu$ writes
\begin{equation} \label{eq:mu intro} -\gamma \Delta_{\mathcal S} \mu(t,s) + \alpha \mu(t,s) = \beta K(t,s)\,,  \end{equation}
where $K  = \kappa_s\kappa_\theta$ is the Gaussian curvature, and $\alpha,\beta,\gamma$ are positive constants.

We assume that the dynamics of release occurs faster than growth, thereby the reaction-diffusion equation \eqref{eq:mu intro} is at  quasi-stationary equilibrium. The choice of the source term $\beta K$ is motivated as follows: 1) the material is released preferentially in the regions of higher curvature, 2) $K$ is an intrinsic invariant of the surface, 3) the whole quantity $\int_{\rm wall}  K$ is constant. Thus, $\mu(t,s)$ can be viewed as the result of redistributing a limited amount of growth material according to the local geometry, together with lateral diffusion on the surface.

Behind this particular choice for the source term,
we have in mind the following more realistic process, involving the distribution of microtubules in fission yeast (see for instance the computational model developed in \cite{Foethke08}). It is well-established that the cytoskeleton controls cell polarity and cell shape \cite{Martin09,terenna_physical_2008}. In normal conditions, microtubules align preferentially in the axis of growth. Thus, they deliver a group of proteins (the +TIP complex) to the cell tip. This enhances the local organization of the growth machinery toward cell tips (see \cite{martin_new_2005} and references therein). Interestingly enough, it is possible to redirect the location of the growth zone by mechanically acting on the cytoskeleton \cite{MincChang09,terenna_physical_2008}. These experiments show evidence of a feedback loop between cytoskeleton organization and the cell shape. Our set of equations can be viewed as a minimal model accounting for this loop.  We ignore secondary feedbacks that regulate the size of the polar cap. We simply assume that microtubules accumulate and deliver the growth material in the region of higher curvature. The growth material is retained in the polarisome. Consequently,  growth is focused at the cell tip. In our basic model, this retention process is associated with the length scale of the diffusion process, namely $(\gamma/\alpha)^{1/2}$.

Previous works studied the maintenance of rod-like cell shape. Generally, the growth pattern is prescribed as a function of the distance to the growing tip of the cell \cite{CampasMahadevan09}, {\it e.g.} via a Gaussian distribution of the growth material along the cell wall. In \cite{DrakeVavylonis13}, the variance of the Gaussian distribution can depend on macroscopic quantities such as the length and the mean radius of the cell. Here, on the contrary, we do not assume that the cell is rod-shaped initially. Moreover, the distribution of the growth material intrinsically depends on the local geometry of the cell wall, {\it e.g.} through its curvature. We focus on the possible initiation of a rod shape from a spherical one, as observed experimentally \cite{Kelly11}. We address mathematically the following morphogenesis question: starting from a small perturbation of a growing radially symmetrical  shape, can the model   evolves towards a rod shape? We fully answer this question using refined linear stability analysis of the radially symmetrical growing shape. 

Other studies have investigated the dynamics of cell growth in prokaryotic actinomycetes. On the contrary to eukaryotic fungi, such as fission yeast, the cytoskeleton plays certainly a minor role in the establishment of the growth pattern: wall building material is likely to be transported to the tip by diffusion. Goriely and Tabor propose two models to investigate self-similar growth of the tip. In \cite{goriely_biomechanical_2003,goriely_self-similar_2003}, they develop a model based on large-deformation elasticity theory. The membrane is described as an axisymmetric elastic shell far from the reference configuration. Similarly to the models discussed above \cite{DrakeVavylonis13,Dumais06}, the inhomogenous elastic modulus is given {\it a priori}, being assumed that the wall gets stiffer far from the tip. In \cite{GorielyGyorgy05}, the same authors develop a purely geometrical model for tip growth, thus neglecting biomechanical effects. The models follows previous studies on the morphogenesis of unicellular algae \cite{pelce_geometrical_1992,pelce_geometrical_1993}. In this model, surface evolution is determined by purely kinematical considerations, and local deposit of material, yielding areal growth. The rate of areal growth is a function of the local geometry ({\it e.g.} the Gaussian curvature). The authors obtain various shapes of self-similar growth in the 2D and in the 3D case. 
\bigskip

%

We can summarize our results as follows. Firstly, we noticed that without taking into account the mechanical effects, the dynamics of wall expansion in the normal direction proportionaly to $\mu$ (or $F(\mu)$) instantaneously creates a very singular pattern. In fact, the system is likely to be ill-posed. Unsurprisingly, taking into account mechanics has a stabilizing effect. We are able to measure quantitatively this effect in the linear regime. The range of parameters for which the system undergoes symmetry breaking is surprisingly narrow. For instance, it is required that $F$ is strongly nonlinear. In addition, $\gamma$, the diffusion coefficient of $\mu$,  has to be relatively small. Then, the morphogenetic instability arises as a competition between the mechanical effects, and the redistribution of material according to the geometry. 

Our analytical work is complemented with numerical simulations far from the radially symmetric shape. 
Interestingly enough, it is unstable, the system selects various anisotropic shapes, including the rod-like shape.

We believe this work paves the way for further (nonlinear) investigations of the coupling between growth and form, and provides some quantitative basis for the analysis of shape selection in walled cells. An interesting perspective would be to adapt this framework to other types of growing walled cells, such as bacteria {\it E. coli} \cite{chang_how_2014}, for which another kind of feedback between cell wall growth and curvature has recently been established \cite{Ursell14}.  

\bigskip

The paper is organized as follows.
Section \ref{Sec:Modeling} deals with 
the complete presentation of the model. 
After introducing some geometrical notations in Subsection \ref{Subsec:Notations}, kinematics are described in Subsection \ref{sub:kinematics} and the mechanical aspects are given in Subsection \ref{sub:mechanics}. The last  Subsection \ref{Subsec:2D_model} introduces an analogous model in the two-dimensional case. It governs the growth of a closed curve, including the main ingredients of the three-dimensional axisymmetric case. 
Section \ref{Sec:Main_Results} presents the main theoretical results of this paper.
We deal with the 2D case in Subsection \ref{subsec:2d}. The linear stability analysis of the three-dimensional case is expressed in Subsection \ref{Subsec:3d}. 
In addition, under the (unrealistic) hypothesis of a fixed length, we characterize in Subsection \ref{subsec:stationnary} the stability of the spherical configurations. 
The proofs are given in Section \ref{sec:Proofs}. 
Section \ref{Sec:Numerics} presents numerical results. 
We confirm the results of the linear stability analysis, and we illustrate the nonlinear dynamics of the system. 

\section{Description of the biomechanical model}
\label{Sec:Modeling}
In this part, we introduce step by step the PDE system which governs the evolution of the cell wall. 
We focus on a model which combines two main conceptual parts.
\begin{enumerate}
\item The cell wall expansion derives from the constitutive laws which govern the deformation of an elastic boundary under pressure, and subject to remodeling. 
We use a standard mechanical model proposed in \cite{CampasMahadevan09,DrakeVavylonis13,Dumais06}. The mechanical parameters are the wall extensibility $\Psi$, the Poisson's ratio $\nu$ and the turgor pressure $P$. 
\item The extensibility depends on the local concentration $\mu$ of some growth material which is distributed along the cell wall. 
The distribution of  material depends intrinsically on the geometry of the cell. We opt for a basic coupling via the Gaussian curvature as in \cite{GorielyGyorgy05}, but the mathematical results obtained here can be generalized to other couplings. 
\end{enumerate}



\subsection{Geometry of the cell wall}
\label{Subsec:Notations}
\begin{figure}[H]
\begin{center}
\includegraphics[width=10cm]{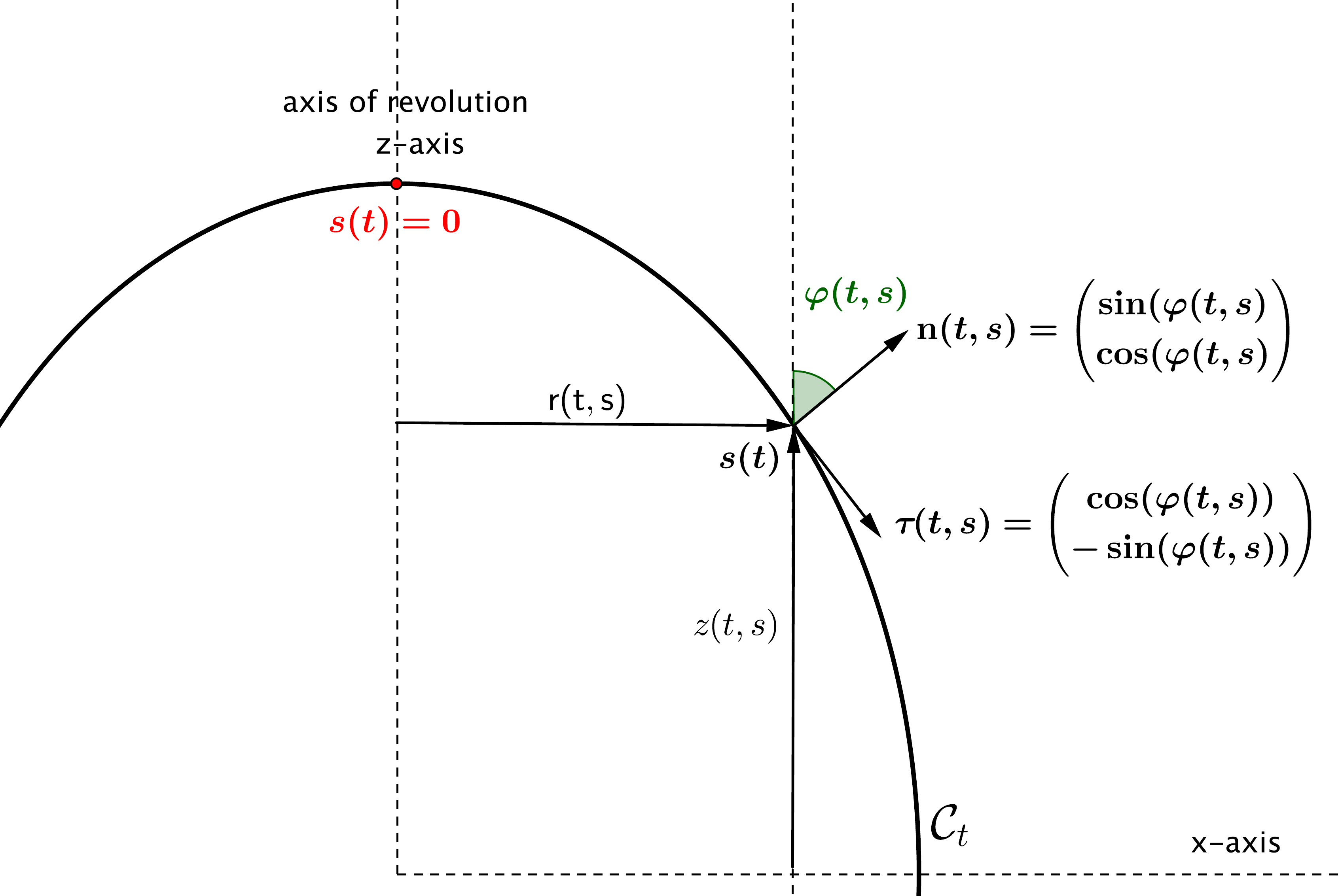}
\caption{\small\label{Fig:surface_revolution} Generatrix curve $\mathcal{C}_t$ of the surface of revolution in the $xz$-plane}
\end{center}
\end{figure}

We restrict to axisymmetric cells. 
The cell wall is a surface of revolution. The generatrix curve in the $x$-$z$ plane is parametrized by $\mathcal{C}_t:=\left\{\left(r(t,s),z(t,s)\right)\right\}$ (see Figure \ref{Fig:surface_revolution}) where $s$ is the curvilinear abscissa. We denote by $L(t)$ the length of the generatrix curve at time $t>0$. 
As described in Figure \ref{Fig:surface_revolution}, $\varphi$ is the angle between the normal  of the curve $\vecn$ and the $z$-axis, while $\vectau$ is the tangent vector. We note $(\vece_1,\vece_2,\vece_3)$ the canonical basis of $\mathbb{R}^3$.
We also denote by $\kappa_s$ and $\kappa_\theta$ the two principal curvatures of the surface.
Since $\mathcal{S}_t$ is a surface of revolution, they have the following expression (see \cite{Carmo76} p. 161),
\begin{equation}
\kappa_s  = \frac{\ds\sqrt{1- (\partial_s r)  ^2}}{r(t,s)}\,,\quad \kappa_{\theta}(t,s) = \frac{-\ds \partial_s r }{\ds\sqrt{1-(\partial_s r)^2}}\,.
\end{equation}


\subsection{Cell wall expansion: kinematics}
\label{sub:kinematics}

For the sake of completeness, we derive the equation governing the evolution of a curve which moves according to a vector field, $\vecv(t,\cdot) : \mathcal{C}_t \to \mathbb{R}^2$.

\begin{figure}[H]
\begin{center}
\includegraphics[scale=0.9]{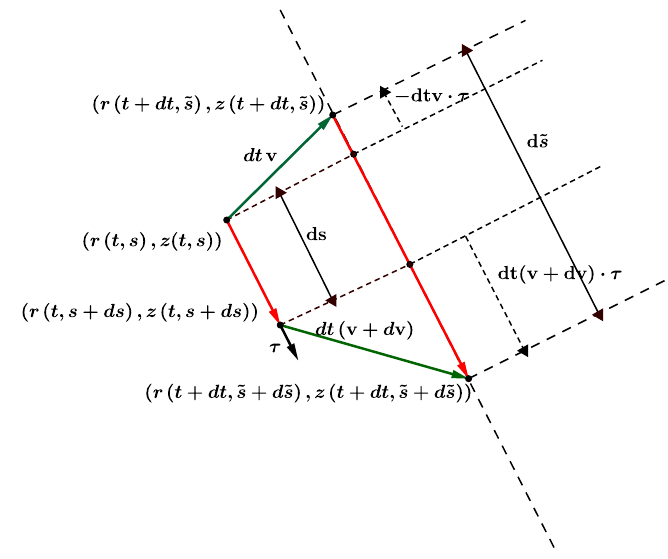}
\caption{\small\label{Fig:eq_evolution} Infinitesimal element, called $\mbox{ds}$, pushed by a vector field $\vecv$.}
\end{center}
\end{figure}

Let $s \in [0,L(t)]$. 
According to Figure \ref{Fig:eq_evolution}, we express the time variation $\frac{\partial }{\partial t}\left(\frac{\partial r}{\partial s}\right)$ as follows, 
\begin{equation*}
\frac{\partial^2 r(t,s) }{\partial t\partial s}   = - \frac{\partial r(t,s)}{\partial s} \left(\frac{\partial \vecv(t,s)}{\partial s} \cdot \vectau(t,s)\right) + \frac{\partial \vecv(t,s)}{\partial s} \cdot \vece_1\,.
\end{equation*}
Moreover, by following the variations around a material point, we get
\begin{equation*}
\frac{\partial }{\partial t}\left(\frac{\partial r(t,s(t))}{\partial s}\right) = \frac{\partial^2 r(t,s(t)) }{\partial t\partial s}    + \frac{ds(t)}{dt} \frac{\partial^2 r(t,s(t))}{\partial s^2}
\end{equation*}
On the other hand, we notice that the curvilinear abscissa is reparamatrized by the elongation of the curve, 
\[
\ds\frac{ds(t)}{dt} = \int_0^{s(t)} \frac{\partial  \vecv(t,s')}{\partial s}  \cdot\vectau(t,s')\, \mbox{ds}'\,,\] 
Finally, the equation governing the evolution of the generatrix curve reads,  for $s\in [0,L(t)]$, as follows 
\begin{multline}
\label{eq:r_final}
\frac{\partial }{\partial t}\left(\frac{\partial r(t,s)}{\partial s}\right) = \ds \frac{\partial}{\partial s} \left( -\left(\ds\int_0^s \frac{\partial \vecv(t,s')}{\partial s}  \cdot\vectau(t,s') \mbox{ds}' \right) \frac{\partial r(t,s)}{\partial s}\right)  
 \ds +\, \frac{\partial  \vecv(t,s)}{\partial s} \cdot\vece_1\,.
\end{multline} 
We refer to Appendix \ref{sec:kinematics} for more details concerning the derivation of \eqref{eq:r_final}.


\subsection{Cell wall expansion: mechanics}
\label{sub:mechanics}


The cell wall dynamics derives from the constitutive laws governing the deformation of an elastic boundary under pressure.
We refer to the biophysical literature for the justification and more in-depth derivation of the equations \cite{CampasMahadevan09, DrakeVavylonis13,Dumais06, Ugural99}.



Following \cite{Dumais06}, we consider the following  constitutive equations relating the strain rates and the stresses,
\begin{equation}
\left\{\begin{array}{ll}
\dot{\epsilon}_s &=  \Psi   \left(\sigma_s -\nu \,\,\sigma_{\theta}\right)\,, \\
\dot{\epsilon}_{\theta} & =\Psi   \left(\sigma_{\theta} -\nu\,\, \sigma_s\right)\,,
\end{array}\right.\label{eq:constituive_laws}
\end{equation}
where
$\sigma_s$ (resp. $\sigma_\theta$) is the meridional stress (resp. the circumferential stress). Parameters are: $\nu$, which is analogous to the Poisson  ratio in linear elasticity, and $\Psi$, the cell wall extensibility. 
We refer to \cite{DrakeVavylonis13} which derive an equivalent model to the one contained in \cite{Dumais06} in a slightly different context (elastic cell wall remodeled under turgor pressure). In \cite{DrakeVavylonis13}, $\Psi$ is linked to the remodeling rate of the cell wall.

The assumption of axial symmetry enables to relate the stresses $\sigma_s$ and $\sigma_{\theta}$ to the principal curvatures $\kappa_s$ and $\kappa_\theta$
\begin{equation}
\left\{\begin{array}{ll}
\ds\sigma_s &= \dfrac{P}{2 \kappa_{\theta}}\,,\smallskip\\
\ds\sigma_{\theta}& =\dfrac{P}{2\kappa_\theta} \left(2 - \dfrac{\kappa_s}{\kappa_\theta}\right)\,. \label{eq:pressure_stress}
\end{array}\right.
\end{equation}
%
%
Kinematic relations for axisymmetric shells make the relation between the strain rates and the velocity.
As usual, 
we decompose the velocity field $\vecv$ as $\vecv =v_n\vecn + v_{\tau} \vectau$. We have \cite{DrakeVavylonis13}
%
\begin{equation}
\left\{\begin{array}{ll}
\ds \dot{\epsilon}_s &= \ds  v_n \kappa_s + \frac{\partial v_{\tau}}{\partial s}\,,\smallskip\\
\ds \dot{\epsilon}_{\theta} &= \ds v_n \kappa_{\theta} + \frac{v_{\tau} \cos(\phi)}{r}\,.\label{eq:strain_velocity}
\end{array}\right.
\end{equation}
Thus, by substituting equation \eqref{eq:constituive_laws} into \eqref{eq:strain_velocity} and by using the relations \eqref{eq:pressure_stress}, we deduce that the functions $v_n$ and $v_{\tau}$ solve the following system of equations on $\mathbb{R}\times [0,L(t)]$,
\begin{subequations}
\label{eq:velocity_stress}
\begin{eqnarray}
\ds \kappa_{\theta} v_n +\frac{v_{\tau} \cos(\phi)}{r} - \Psi (\sigma_{\theta} -\nu\sigma_s)=0\,,\\ \label{eq:vn} 
\ds \frac{\partial v_{\tau}}{\partial s} - \left( \frac{\kappa_s}{\kappa_{\theta}}\frac{\cos(\phi)}{r}\right) v_{\tau} + \frac{\kappa_s}{\kappa_{\theta}}\Psi\left(\sigma_{\theta}-\nu\sigma_{s}\right)-\Psi \left(\sigma_s-\nu\sigma_{\theta}\right)=0\,.
 \label{eq:vt} 
\end{eqnarray}
\end{subequations}

\subsection{Cell wall expansion: growth pattern}

We assume that the function $\Psi$, which can be read as the cell wall extensibility, or the remodeling rate, depends on some growth factor released  in the vicinity of the cell wall. The surface concentration of this growth factor is denoted by $\mu(t,s)$. We assume a general relationship between $\Psi$ and $\mu$, 
\[\Psi(t,s)=F(\mu(t,s))\,, \] 
where $F$ is an increasing function.
We suppose that the growth factor $\mu$ is diffusing laterally, and is released locally, proportionally to  the Gaussian curvature. Moreover, we assume that the dynamics of release-diffusion are much faster than the time scale of remodeling and growth. Thus, the density $\mu$ satisfies the following equation
\begin{equation}
\label{eq:mu}
-\gamma\,\Delta_{\mathcal{S}} \mu(t,s) +\alpha\, \mu(t,s) =  \beta   \kappa_s(t,s) \kappa_\theta(t,s) \,,
\end{equation}
where $\gamma$,$\alpha$ and $\beta$ are positive constants and $\Delta_{\mathcal{S}}$ is the Laplacian-Beltrami operator on the cell wall ${\mathcal{S}}$.
In the case of an axisymmetric surface, the Laplacian-Beltrami operator $\Delta_{\mathcal{S}}$ of $\mu$  is expressed as follows \cite[p. $79$]{Carmo76}. 
\begin{equation}
\label{eq:Laplacian_Beltrami}
- \Delta_{\mathcal{S}} \mu(t,s) = - \frac{1}{r(t,s)}\left(\dfrac{\partial}{\partial s} \left(r(t,s)\,\dfrac{\partial \mu (t,s)}{\partial s}  \right)\right)\,.
\end{equation}

At this stage, we can compute the rate of expansion of a spherical shape. 
  
\begin{proposition}[Spherical solution]
\label{Prop:3D_spheres}
Let $L^c$ be the solution of the following ODE
\begin{equation}
\label{eq:3D_L_c}
\dfrac{dL^c(t)}{dt}  =\frac{\pi}{2} P \left( 1-\nu \right)F \left( \frac{\beta\pi^2}{\alpha  L^c(t)^2}  \right) L^c(t)^2\,, 
\end{equation}
with a given initial length $L^c(0)$.
Then, $\left(r^c,\mu^c\right)$ defined by
\begin{equation}
\left\{\begin{array}{ll}
r^c(t,s) &= \dfrac{L^c(t)}{\pi} \sin\left(\dfrac{s\pi}{L^c(t)}\right)\,,\medskip\\
\mu^c(t) & =  \ds\frac{\beta\pi^2}{\alpha  L^c(t)^2}\,\label{eq:circle}\, ,
\end{array}\right.
\end{equation}
is a particular solution of the system \eqref{eq:r_final}-\eqref{eq:velocity_stress}-\eqref{eq:Laplacian_Beltrami}. 
 \end{proposition}


%
%
%
%

In Section \ref{Sec:Main_Results}, we address the issue of stability of the solution \eqref{eq:circle}.


\subsection{A simplified 2D model}
\label{Subsec:2D_model}

In this subsection, we derive the dynamics which governs a  thin elastic string under pressure. More precisely, we assume that the cell wall is  a closed curve which is symmetric with respect to the $z$-axis. 

We first express the Hooke's law,
$\varepsilon = \Psi \sigma$, 
where $\varepsilon$ is the strain rate and $\sigma$ is the stress acting on the curve. The latter is determined by  the Laplace  law, $\sigma = {P}/{\kappa}$, where $\kappa$ is the curvature.
Moreover, by geometrical consideration, the strain rate is also equal to
$\varepsilon = v_n \kappa$,
where $v_n$ is the normal component of the velocity vector field.
All in all, the normal velocity of the curve is given by,
\begin{equation}
\label{eq:strain3}
v_n=\Psi\, \frac{P}{\kappa^2}\,.
\end{equation}
We use the same notations as for the 3D axysimmetric case. In particular, we have $\kappa = -\partial_s \phi$. By analogy with the 3D case,  the evolution of the curve obeys the following equation
\begin{multline}
\frac{\partial }{\partial t}\left(\frac{\partial r(t,s)}{\partial s}\right)=- \frac{\partial }{\partial s} \bigg(\left(\ds\int_0^s P\, \frac{\Psi(t,u)}{\partial_s\phi(t,s')}\,\mbox{ds}' \right) \frac{\partial r(t,s)}{\partial s}\bigg) \\
\,+\frac{\partial }{\partial s}  \bigg( \ds P \frac{\Psi(t,s)}{\left(\partial_s\phi(t,s)\right)^2} \sin(\phi(t,s))\bigg)\,.\label{eq:2D_r}
\end{multline}
The coefficient $\Psi$ is supposed to depend on some  growth material which is released proportionally to the curvature, and diffuses quickly along the curve,
\begin{equation}
-\gamma\frac{\partial^2 \mu(t,s)}{\partial s^2}  +\alpha \mu(t,s) = \beta \frac{\partial \phi(t,s)}{\partial s}  \,.\label{eq:2D_mu}
\end{equation}


The next proposition states that the circular shape is preserved by the system \eqref{eq:2D_r}-\eqref{eq:2D_mu}.

\begin{proposition}[Circular solution]
\label{Prop:2D_spheres}
Let  $L^c$ be the solution of the following ODE
\begin{equation}
\label{eq:2D_circle_lenght}
\dfrac{dL^c(t)}{dt}  = \ds \frac{P}{\pi}F\left(\frac{\beta\pi}{\alpha  L^c(t)}\right)L^c(t)^2\,
\end{equation}
with a given initial length $L^c(0)$.
Then, $\left(r^c,\mu^c\right)$ defined by
\begin{equation*}
\left\{\begin{array}{ll}
r^c(t,s) &= \ds \frac{L^c(t)}{\pi} \sin\left(\frac{s\pi}{L^c(t)}\right)\,,\medskip
\\
\mu^c(t) &= \ds\frac{\beta\pi}{\alpha  L^c(t)}\,, \end{array}\right.
\end{equation*}
 is a solution to the system of equations \eqref{eq:2D_r}-\eqref{eq:2D_mu}. 
 \end{proposition}



\section{Linear stability analysis}
\label{Sec:Main_Results}


In this section, we derive the precise conditions under which the cell wall expansion is unstable around the radially symmetric solution. The result is given  in Subsection \ref{subsec:2d}  for the  $2$D case and in Subsection \ref{Subsec:3d} for the $3$D case. Note that the 2D model cannot be viewed as a particular case of the 3D system, but there are strong similarities between the two settings. 

To analyze the stability of the model, we proceed in two steps. The first step consists in transforming the original system into a more tractable one by using an appropriate change of variables and by making the physical quantities dimensionless.
The second step is to solve the linearized equation around the circular (2D case) or the spherical (3D case) solution.


The following definition will be useful.
\begin{definition}[Degree of nonlinearity]
\label{def:degree}
\noindent For $f:\R_+\to \R$ a differentiable function, and $\mu\in \R_+$, we define the degree of  nonlinearity $\deg(F; \mu)$ of $f$ at point $\mu$ as following
$$
\deg(F;\mu)= \frac{f'(\mu) \mu}{f(\mu)} = \dfrac{d \log f(\mu)}{d\log(\mu)}\,.
$$
\end{definition}
It coincides with the notion of elasticity in economics. For a homogeneous function $f(\mu) = \mu^p$, we have $\deg(F;\mu)=p$.

In the sequel, the curve is parametrized by the angle $\phi$. The following relation enables to rewrite the equations for the wall expansion \eqref{eq:r_final} and \eqref{eq:2D_r},
\begin{equation}
\label{Change_Variable}
\frac{\partial r (t,s)}{\partial s} =  \cos(\phi(t,s))\,.
\end{equation}

\subsection{Stability results for the 2D model}
\label{subsec:2d}

Firstly, we adimensionalize the equations. Let $\tilde{\varphi}$ and $\tilde{\mu}$ be defined as 
$$
\tilde{\varphi}(t,x) = \varphi(t,xL(t))  \quad \textrm{and}\quad  \tilde{\mu}(t,x)= \dfrac{\alpha L(t)}{\beta \pi}\mu(t,xL(t)) = \dfrac{L(t)}{L^c(t)}\dfrac{\mu(t,xL(t))}{\mu^c(t)}\,, \quad t>0\, , \quad x\in (0,1)\, .
$$
We introduce $\tilde{\Psi}(t,x)=  F\left(\frac{L^c(t)}{L(t)}\mu^c(t)\tilde\mu(t,x)\right)$. 
We also introduce the reduced diffusion coefficient on the membrane,
\[\sigma(t)=\frac{\gamma \pi^2}{\alpha L(t)^2}\,. \]
Then, 
$\left(\tilde{\varphi} ,\tilde{\mu} \right)$ is solution to the following system of equations
\begin{multline}
\label{eq:2D_new_Phi}
\partial_t(\tilde{\varphi}(t,x)) =  \left(\ds\frac{L'(t)}{L(t)}x - PL(t) \ds \int_0^x \frac{\tilde{\Psi}(t,u)}{\partial_x\tilde{\varphi}(t,u)}\mbox{du}\right)\partial_x\tilde{\varphi}(t,x)\\
-PL(t) \left(\ds\frac{\partial_x\tilde{\Psi}(t,x)}{\left(\partial_x\tilde{\varphi}(t,x)\right)^2}-  \frac{2\tilde{\Psi}(t,x)}{\left(\partial_x\tilde{\varphi}(t,x)\right)^3}\partial_x^2\tilde{\varphi}(t,x)\right)\,,
\end{multline}
\begin{equation}
\label{eq:2D_new_mu}
-\frac{\sigma(t)}{\pi^2}\partial_x^2 \tilde{\mu}(t,x) +   \tilde{\mu}(t,x)=  \frac1\pi \partial_x \tilde{\varphi}(t,x) \,,
\end{equation}
with the boundary conditions $\partial_x \tilde{\mu}(t,0)=\partial_x \tilde{\mu}(t,1) = 0$. 
%
%
The following result quantifies the local stability of the system \eqref{eq:2D_new_Phi}-\eqref{eq:2D_new_mu} around the circular solution. 

\begin{theorem}[Linear stability analysis]
\label{thm:2D}
Let $(u,v)$ be the solution of the linearized version of system \eqref{eq:2D_new_Phi}-\eqref{eq:2D_new_mu} around the circular solution. Let $(a_k)_{k\geq 1}$ denote the Fourier coefficients of $\partial_x u$, defined as 
$
\partial_x u(t,x) = \ds \sum_{k\geq1}  a_k(t) \cos(k\pi x)
$.
The coefficients are solutions to the following set of linear ODE 
$$
a'_k(t) = M(t,k) a_k(t)\,,
$$
where $M$ is a rational fraction of  $k$, the sign of which is determined by the sign of the following polynomial for any integer value $k\geq 1$,
$$
N(t,k) = \Big[1 - \deg(F; \mu^c(t))\Big] + \Big[\sigma^c(t) - 2    + \deg(F;\mu^c(t))\Big] k^2 -  \Big[2 \sigma^c(t)\Big]  k^4  \,.
$$ 
\end{theorem}
The proof of this result is given in Subsection \ref{2D_Proof}. We define the functions $G_1$ and $G_2$ as,
\begin{align*}
G_1(\sigma,d) &= 1+\left(d-2\right)\frac{1}{\sigma}\\
G_2(\sigma,d) &= 1+ 2\left(2-3d\right)\frac{1}{\sigma} + \left(d-2\right)^2 \frac{1}{\sigma^2}\,.
\end{align*}
For a given $t>0$, the zeros of $N(t,\cdot)$ are expressed as 
$$
\pm \frac{1}{2} \sqrt{G_1\left(\sigma^c(t), \deg(F;\mu^c(t))\right)\pm \sqrt{G_2\left(\sigma^c(t), \deg(F;\mu^c(t))\right)}}\,.
$$
The consequences of this result are discussed in Section \ref{subsec:stationnary} and illustrated in Section \ref{Sec:Numerics}.


\subsection{Stability results for the 3D model}
\label{Subsec:3d}

In this Subsection, the previous results are extended to the 3D case. We proceed with the same nondimensionalization.
%
Let $\tilde{\varphi}$ and $\tilde{\mu}$ be defined as 
$$
\tilde{\varphi}(t,x) = \varphi(t,xL(t))  \quad \textrm{and}\quad  \tilde{\mu}(t,x)= \dfrac{\alpha L(t)^2}{\beta \pi^2} \mu(t,xL(t))  =  \dfrac{L(t)^2}{L^c(t)^2}\dfrac{\mu(t,xL(t))}{\mu^c(t)}\,, \quad t>0\, , \quad x\in (0,1)\, .
$$
We introduce $\tilde{\Psi}(t,x)=  F\left(\frac{L^c(t)^2}{L(t)^2}\mu^c(t)\tilde\mu\right)$. 
We also introduce the reduced diffusion coefficient on the membrane,
\[\sigma(t)=\frac{\gamma \pi^2}{\alpha L(t)^2}\,, \]
and the auxiliary functions
$\mbox{Icos}$ and $\mbox{Icosin}$, defined as
\begin{equation}
\ds\mbox{Icos}(t,x)=\int_0^x \cos(\tilde{\varphi}(t,x'))\mbox{dx}'  \quad \textrm{and}\quad  \mbox{Icosin}(t,x)=\frac{\mbox{Icos}(t,x)}{\sin(\tilde{\varphi}(t,x))}\,.
\label{eq:Icos}
\end{equation}
Then, $\left(\tilde{\varphi}(t,x),\tilde{\mu}(t,x)\right)$ is solution to the system of equations
\begin{multline}
\label{eq:3D_new_Phi}
\partial_t \tilde{\varphi}(t,x)= \ds \left( \ds\frac{L'(t)}{L(t)}  x
- \int_0^x \left(d(t,x') +  H(t,x')  \partial_x\tilde{\varphi}(t,x')\right)\mbox{dx}'\right)\,\,\partial_x \tilde{\varphi}(t,x)  \\
 + \cot(\tilde{\varphi}(t,x)) \,d(t,x) -\partial_x H(t,x)\,,
\end{multline}
 \begin{align}
\label{eq:3D_new_mu}
  - \dfrac{\sigma(t)}{\pi^2} \left( \partial^2_x \tilde{\mu}  + \frac{\cos(\tilde{\varphi})}{\mbox{Icos}}\partial_x\tilde{\mu} \right) + \tilde{\mu}  = \frac1{\pi^2}  \partial_x\tilde{\varphi}\,   \frac{\sin(\tilde{\varphi})}{\mbox{Icos}}  \,,
\end{align}
where the function $d$ and $H$ are given by 
\begin{align}
\Lambda_1 &=\frac{1}{2} P L  \tilde{\Psi} \left( (1-2\nu)  \mbox{Icosin} + 2\left(\nu-1\right) \left(\mbox{Icosin}\right)^2\partial_x \tilde{\varphi} 
 +  \left(\mbox{Icosin}\right)^3 \left(\partial_x\tilde{\varphi}\right)^2\right)\,,\label{eq:def d}\\
 \Lambda_2 &= \frac{1}{2} P L \tilde{\Psi} \left( \left(2-\nu \right) \left(\mbox{Icosin}\right)^2 -\left(\mbox{Icosin}\right)^3\partial_x \tilde{\varphi} \right)\label{eq:def H}\,.
\end{align}

The details of the change of variables are given in Appendix \ref{Appendix:Lemma}. 
The following result quantifies the local stability of the system  \eqref{eq:3D_new_Phi}-\eqref{eq:3D_new_mu} around the spherical solution.


\begin{theorem}[Linear stability analysis]
\label{thm:3D}
Let $(u,v)$ be the solution of the linearized version of system \eqref{eq:3D_new_Phi}-\eqref{eq:3D_new_mu}around the spherical solution. Let $(a_k)_{k\geq 1}$ denote the Fourier coefficients of $\partial_x u$, defined as $
\partial_x u(t,x) = \ds \sum_{k\geq1}  a_k(t) \cos(k\pi x)$.
The subset of coefficients $A(t) = (a_k(t))_{k\geq 2}$ are solutions to a triangular system of ODE, denoted by 
\begin{equation}
\label{linearized}
A'(t) = M(t) A(t)\,, 
\end{equation}
where the signs of the diagonal coefficients of $M$ are determined by   the sign of the following polynomial for any integer value $k\geq 2$,
\begin{equation}\label{eigenval}
N(t,k) = n_0(t) + n_1(t)  k + n_2(t)  k^2 + n_3(t)  k^3 +n_4(t)  k^4 \,,\quad \forall t>0\,, 
\end{equation}
with
\begin{align*}
n_0(t) &=  (1-\nu )\left[1-2 \deg(F;\mu^c(t))\right] \,,\\
n_1(t)  &= \,-(1-\nu)\left[\sigma^c(t) + \deg(F;\mu^c(t))\right] +1 \,,\\
n_2(t)  &=  -\nu\sigma^c(t) -1 + \left( 1-\nu \right) \deg(F;\mu^c(t)) \,,\\
n_3(t) &= 2\sigma^c(t)  \,,\\
n_4(t) &= -\sigma^c(t) \,.
\end{align*}
\end{theorem}

The proof of this result is given in Subsection \ref{3D_Proof}. We define the functions $G_1$ and $G_2$ as
\begin{align*}
G_1(\sigma,d) &= 3-2\nu + 2\left(\left(1-\nu\right)d - 1\right)\frac{1}{\sigma}\,,\\
G_2(\sigma,d) &= (1-\nu)^2 -2(1-\nu)((3+\nu)d-1)\frac{1}{\sigma} + \left((1-\nu)d - 1\right)^2  \frac{1}{\sigma^2}\,.
\end{align*}
For a given $t>0$, the zeros of $N(t,\cdot)$ are expressed as 
$$
\frac{1}{2} \pm \frac{1}{2}\sqrt{G_1(\sigma^c(t), \deg\left(F;\mu^c(t)\right)) \pm  2\sqrt{G_2(\sigma^c(t), \deg\left(F;\mu^c(t)\right))}}\, .
$$

\subsection{Instability in the case of quasi-stationnary expansion}
\label{subsec:stationnary}

In this Subsection, we discuss several necessary conditions for instability based on the previous analysis. 

We make two assumptions to simplify the discussion. Firstly, we do not restrict $k$ to taking integer values. We assume that it can be any real value greater than unity. Secondly, we make the modelling assumption that the dynamics of (in)stability is faster than wall expansion. Therefore, we consider that $L(t)$ varies slowly, so slowly that we further assume that it takes a constant value $L$.  

\paragraph{The 2D case.} We begin with the two-dimensional case. There are four possible roots, which can be real or not by pairs. In addition, we require that the largest root is greater than one. 

We find that the latter condition implies $G_1 + \sqrt{G_2} \geq 4 $, which is satisfied if $G_1\geq 4$ or if $G_2\geq (4 - G_1)^2$. The latter cannot be satisfied since $G_2 - (4 - G_1)^2 = -8(1 + 1/\sigma)<0$, therefore we obtain the necessary condition $G_1\geq 4$, which reads
\[ d\geq 2 + 3\sigma\, . \]
In particular, it is mandatory that the degree of nonlinearity is larger than 2.  This could have been guessed from the expression of the expansion velocity 
\eqref{eq:strain3}. In fact, if $\Psi = \kappa^2$, then the velocity of expansion is constant, and the geometric coupling disappears. 

On the other hand, a simple asymptotic analysis of the roots in the two cases i) $d\to +\infty$, and ii) $\sigma\to 0$  together with the conditon $(d>2)$, shows that high degree of nonlinearity and small lateral diffusion on the cell wall guarantee instability. In fact, by introducing the notation $k_1 = \frac12\sqrt{G_1 + \sqrt{G_2}}$ and $k_2 = \frac12\sqrt{G_1 - \sqrt{G_2}}$, we notice that i) as $d\to +\infty$, we have respectively $k_1 \to +\infty $ and $k_2 \to 1$, and ii) as $\sigma\to 0$, we have respectively $k_1 \to +\infty $ and $k_2 \to \sqrt{1 + \frac1{d-2}}>1 $.

\paragraph{The 3D case.} We perform the same analysis in the three-dimensional case. Notations are the same as in the 2D case. We find that the condition $k_1\geq 2$ implies $G_1 + 2\sqrt{G_2} \geq 9 $, which is satisfied if $G_1\geq 9$ or if $4 G_2\geq (9 - G_1)^2$. The latter cannot be satisfied   since $4G_2 - (9 - G_1)^2 = -16(1  + \nu)(2 + 1/\sigma)<0$. Therefore, we obtain the necessary condition $G_1\geq 9$, which reads also
\begin{equation}   
d \geq \dfrac{1 + (3+\nu)\sigma}{1 - \nu}\, . 
\label{eq:necessary}
\end{equation}

The same asymptotic analysis of the roots holds true, provided $(1-\nu)d - 1>0$, which is compatible with the above necessary condition. In particular, i) as $d\to +\infty$, we have respectively $k_1 \to +\infty $ and $k_2 \to 2$, and ii) as $\sigma\to 0$, we have respectively $k_1 \to +\infty $ and $k_2 \to \frac1 2+\frac1 2\sqrt{ 9 + 4\frac{1 + \nu}{(1-\nu)d - 1}}>2$.

\section{Proofs of the linear stability results}
\label{sec:Proofs}


This Section is devoted to the proofs of  Theorem \ref{thm:2D} and Theorem  \ref{thm:3D}, respectively. 
Although the three-dimensional case requires extra computational tricks, the proofs are sketchily the same. We drop the $\;\tilde{}\;$  in the adimensional system for the sake of notation. 

%
%
%
%
%

In the following, we denote respectively by $\varphi^w$, $\mu^w$, and $L^w$ small perturbations of the radially symmetric shape associated whith the particular solution $\varphi(t,x) = \pi x$, $\mu(t,x) = 1$, and $L^c(t)$.


We use the Fourier basis associated with periodic functions of period two, up to extending the functions for $x\in (0,2)$ using axisymmetry, {\it e.g.} 
\begin{equation}
\label{eq:hyp_symmetry}
\varphi(x) = 2\pi - \varphi(2-x)\quad \textrm{and} \quad \mu(x) = \mu(2-x)\,, \;\forall x\in (0,2)\,.
\end{equation}

\subsection{Proof of Theorem \ref{thm:2D} : $2$D case}
\label{2D_Proof}

We aim at linearizing the system  \eqref{eq:2D_new_Phi}-\eqref{eq:2D_new_mu} around the circular shape solution $(\varphi,\mu) = (\pi x, 1)$. 
%
%
The linearized equation of \eqref{eq:2D_new_Phi} is 
\begin{multline}
\label{eq:2D_linearized_1}
\partial_t\phi^w  = - PL^c(t) \left( \ds \int_0^x\Psi^w \mbox{dx}' + \frac{1}{\pi^2}\partial_x\Psi^w \right)   \\
+ PL^c(t)\Psi^c(t)   \left( - \frac{x}{\pi} \partial_x\phi^w  + \ds\int_0^x \frac{1}{\pi}\partial_x\varphi^w  \mbox{dx}'+\frac{2}{\pi^3}\partial_x^2\varphi^w \right) \\
+\left(\frac{L'(t)}{L(t)}\right)^c \,x\, \partial_x \varphi^w  + x\, \pi\, \left(\frac{L'(t)}{L(t)}\right)^w\,.
\end{multline}
Since $L^c$ solves equation \eqref{eq:2D_circle_lenght}, 
equation \eqref{eq:2D_linearized_1} becomes 
\begin{multline}
\label{eq:2D_linearized_1}
\partial_t\phi^w  = - PL^c(t) \left( \ds \int_0^x\Psi^w \mbox{dx}' + \frac{1}{\pi^2}\partial_x\Psi^w \right)   \\
+ PL^c(t)\Psi^c(t)   \left(   \ds\int_0^x \frac{1}{\pi}\partial_x\varphi^w  \mbox{dx}'+\frac{2}{\pi^3}\partial_x^2\varphi^w \right)  
 + x\, \pi\, \left(\frac{L'(t)}{L(t)}\right)^w\,.
\end{multline}
%
From the definition of $\Psi(t,x)$ in the nondimensional variables (see Section \ref{subsec:2d}), we have
\begin{equation} 
\label{eq:diff Psi}
\Psi^w(t,x) =  \mu^c(t)F'(\mu^c(t)) \left( \mu^w(t,x) - \dfrac{L^w(t)}{L^c(t)}\right) =  \deg\left(F;\mu^c(t)\right) F(\mu^c(t)) \left( \mu^w(t,x) - \dfrac{L^w(t)}{L^c(t)}\right)\, .  \end{equation}
By differentiating \eqref{eq:2D_linearized_1} with respect to $x$, we get
\begin{align}
 \partial_t\left(\partial_x\phi^w(t,x)\right) = & - PL^c(t) \deg(F; \mu^c(t)) F(\mu^c(t)) \left( \ds \mu^w(t,x) + \frac{1}{\pi^2}\partial^2_x\mu^w(t,x)\right) \nonumber \\
& + PL^c(t)  F(\mu^c(t)) \left( \ds\frac{1}{\pi}\partial_x\varphi^w(t,x)+\frac{2}{\pi^3}\partial_x^3\varphi^w(t,x)\right) \nonumber\\
  & +  \pi \left(\frac{L'(t)}{L(t)}\right)^w   +  P  \deg(F; \mu^c(t)) F(\mu^c(t))  L^w(t) \,.
  \label{eq:2D_linearized_Phi}
\end{align}
On the other hand, equation \eqref{eq:2D_new_mu} is already linear, so we write
\begin{equation}
\label{eq:2D_linearized_mu}
 - \frac{\diff(t)}{\pi^2} \partial_x^2 \mu^w(t,x) +    \mu^w(t,x)    = \frac1\pi \partial_x \varphi^w(t,x) \,.\\
\end{equation}
We  expand the $2$-periodic functions $\partial_x \varphi^w$ and $\mu^w$  in Fourier series.
Using the symmetric relations   \eqref{eq:hyp_symmetry}, we can write
\begin{align*}
\label{eq:2D_Fourier}
\partial_x\varphi^w(t,x) = \ds  \sum_{k\geq1} a_k(t) \cos(k \pi x)\,,\\
\mu^w(t,x) = \ds \sum_{k\geq1} c_k(t) \cos(k \pi x)\,,
\end{align*}
where 
\begin{align*}
a_k(t) &= \ds 2\int_0^1  \cos(k\pi  x) \partial_x\varphi^w(t,x) \mbox{dx} \,,\\
c_k(t) &= \ds  2\int_0^1  \cos(k\pi  x)\mu^w(t,x) \mbox{dx}\,.
\end{align*}
By replacing the functions $\partial_x\varphi^w$ and $\mu^w$ by their Fourier series in equation \eqref{eq:2D_linearized_mu} we get
\begin{equation}
\label{eq:2D_Fourier_RD}
c_k(t) = \frac{1}{\pi\left(1 + \sigma^c(t) k^2 \right)} a_k(t)\,, \quad \forall k\geq 1.
\end{equation}
Similarly, for equation \eqref{eq:2D_linearized_Phi}, we obtain,
\begin{equation*}
\label{eq:2D_Fourier_Meca}
\dfrac {d a_k}{dt}   = \frac{1}{\pi} P L^c F \left( - \pi \deg(F;\mu^c) \left(1- k^2\right)c_k  +    \left(1 - 2k^2\right) a_k \right) \,.
\end{equation*}
Therefore, we deduce that for all $k\geq1$, the following differential equation holds,
\begin{equation*}
\label{eq:2D_Fourier}
\dfrac {d a_k}{dt}  = \dfrac{  P L^c F}{\pi \left(1 + \sigma^c  k^2\right)} \left( (1 - \deg(F;\mu^c)) + \left( \sigma^c  - 2    + \deg(F;\mu^c)\right)k^2 - 2 \sigma^c  k^4 \right) a_k  \,,
\end{equation*}
which concludes the proof.

%
%
%
%

\subsection{Proof of Theorem \ref{thm:3D} : $3$D case}
\label{3D_Proof}

We aim at linearizing the system  \eqref{eq:3D_new_Phi}-\eqref{eq:3D_new_mu} around the spherical  shape solution $(\varphi,\mu) = (\pi x, 1)$. 
%
%
%
\medskip

\noindent\textbf{Step \#1 : Linearized equations.}
We start by expressing the expansion of functions $d$ \eqref{eq:def d} and $H$ \eqref{eq:def H} around the spherical solution. The zeroth-order term $H^c$ reads
\begin{align*}
\Lambda_2^c =& \left( 1-\nu \right)\frac{P}{2\pi^2} L^c  \Psi^c .
\end{align*}
The first-order term $H^w$ is
\begin{multline*}
       \Lambda_2^w =   \left( 1-\nu \right) \frac{P}{2\pi^2} \left[L^c  \Psi^w+  \Psi^c L^w \right]\\
         + \frac{P}{2\pi^2} L^c \Psi^c\left[\frac{\left(2\nu-1\right)}{\sin(\phi^c)}\int _{0}^{x}\!\cos \left(\varphi^c\right) \partial_x\varphi^w \mbox{dx} - \frac{1}{\pi}\,\partial_x \varphi^w\right]\,.
\end{multline*}
Similarly, the expansion of $d$ writes
\begin{align*}
        \Lambda_1^c =&   0\,,\\
        \Lambda_1^w =&  \nu  \frac{P}{\pi} L^c \Psi^c \, \left[ -\frac{1}{\sin \left( \varphi^c \right) }\int _{0}^{x}\!\cos \left(\varphi^c\right) \partial_x\varphi^w \mbox{dx} + \frac{1}{\pi}\partial_x \varphi^w \right]\,.
 \end{align*}
Then, after differentiating  the linearized version of equation \eqref{eq:3D_new_Phi} with respect to $x$, we get
 \begin{align}
 \label{eq:linearized_1}
 \partial_t \left(\partial_x \varphi^w\right) =& \left( \left(\frac{ L'}{L}\right)^w - \Lambda_1^w - \Lambda_2^c \partial_x\varphi^w - \Lambda_2^w \pi\right) \pi \nonumber\\
 & \quad - \frac{\pi}{\sin(\varphi^c)^2}   \Lambda_1^w + \cot(\varphi^c) \partial_x \Lambda_1^w - \partial^2_x \Lambda_2^w\,.
 \end{align}
The equivalent of \eqref{eq:diff Psi} in the 3D case reads 
\begin{equation} 
\label{eq:diff Psi 3D}
\Psi^w(t,x) =    \deg\left(F;\mu^c(t)\right) F(\mu^c(t)) \left( \mu^w(t,x) - 2 \dfrac{L^w(t)}{L^c(t)}\right)\, .  \end{equation}
We obtain,
\begin{align}
\partial_t \left(\partial_x \varphi^w\right) = & PL^cF(\mu^c)   \left( \frac{1}{\sin(\varphi^c)^3} \ds \int _{0}^{x} \cos \left(\varphi^c\right) \partial_x \varphi^w \,\mbox{dx} \right. \nonumber \\
& + \frac{1}{2\pi}\left(1-\nu - \frac{2}{\sin^2(\varphi^c)}\right)\partial_x \varphi^w 
\left. + \frac{\cos(\varphi^c)}{2\pi^2 \sin(\varphi^c)}\partial^2_x \varphi^w +  \frac{1}{2\pi^3} \partial^3_x \varphi^w\right) \nonumber\\
& - PL^c\deg(F;\mu^c) F   \left(\frac{1}{2} \left(1-\nu\right)\mu^w+\frac{1}{2\pi^2} \left(1-\nu\right)\partial_x^2\mu^w\right) \nonumber \\
& + \pi  \left(\frac{L'}{L}\right)^w  + \frac{1}{2}\left(1-\nu\right)P F \left(2 \deg(F;\mu^c) - 1  \right)L^w \, . \label{eq:linearized_1_bis}
 \end{align}
Then, we multiply both side of equation \eqref{eq:linearized_1_bis} by $\sin(\varphi^c)^3$ and we differentiate  with respect to $x$. Eventually, we divide the resulting equation on both side by $\sin(\varphi^c)^2$, so as to  obtain
\begin{align}
&\partial_t \left(3\pi \cos(\varphi^c) \partial_x \varphi^w+ \sin(\varphi^c) \partial^2_x \varphi^w\right) \nonumber\\
& =   P L^c F   \Bigg(\frac{3}{2} (1-\nu) \cos(\varphi^c)\, \partial_x \varphi^w -\frac{1}{2\pi} (2+\nu)\sin(\varphi^c)  \partial^2_x \varphi^w\nonumber\\
&  \qquad\qquad\qquad\qquad\qquad\qquad +   \frac{2}{\pi^2}  \cos(\varphi^c) \partial^3_x \varphi^w + \frac{1}{2\pi^3} \sin(\varphi^c)\partial^4_x \varphi^w\Bigg) \nonumber\\
& \quad - P L^c \deg(F;\mu^c) F   (1-\nu) \Bigg(\frac{3}{2}\pi \cos(\varphi^c)  \mu^w + \frac{1}{2} \sin(\varphi^c)  \partial_x \mu^w \nonumber\\
&  \qquad\qquad\qquad\qquad\qquad\qquad  + \frac{3}{2\pi} \cos(\varphi^c)  \partial^2_x \mu^w + \frac{1}{2\pi^2} \sin(\varphi^c) \partial^3_x \mu^w\Bigg) \nonumber\\
& \quad +\frac{3}{2} (1-\nu) P F \left(2 \deg(F;\mu^c) - 1  \right)\pi\cos(\phi^c) L^w  + 3 \pi^2 \cos(\phi^c) \left(\frac{L'}{L}\right)^w \,. \label{eq:linearized_2}
 \end{align}
By multiplying the linearized version of equation \eqref{eq:3D_new_mu} by $\sin(\varphi^c)$ and by differentiating both side with respect to $x$, 
we get  
\begin{multline}
\label{eq:mu_linearize}
\lefteqn{2  \cos(\varphi^c) \partial_x \varphi^w +  \frac1\pi \sin(\varphi^c) \partial_x^2 \varphi^w}  \\
 =  \pi \cos(\varphi^c)\mu^w + (\sigma^c + 1) \sin(\varphi^c) \partial_x \mu^w    -2 \frac{\sigma^c}\pi \cos(\varphi^c) \partial_x^2 \mu^w - \frac{\sigma^c}{\pi^2}\sin(\varphi^c) \partial_x^3 \mu^w \,.
 \end{multline}
%

\noindent \textbf{Step \#2: Fourier series expansion.}
%
We  expand the $2$-periodic functions $\partial_x \varphi^w$ and $\mu^w$  in Fourier series. Using the symmetric relations  \eqref{eq:hyp_symmetry}, the Fourier series read as in the 2D case. 
We will need the following formulas: 
\begin{multline}
\label{eq:fourrier1}
\cos(\varphi^c) \partial_x\varphi^w =  a_1 + \left(a_0+a_2\right) \cos(\pi x)  \\
 + \,\,\ds \sum_{k\geq2} \left(a_{k+1}+a_{k-1}\right) \cos(k\pi x)\,. 
\end{multline}
\begin{multline}
\label{eq:fourrier2}
\sin(\varphi^c) \partial^2_x\varphi^w = - \pi a_1 - 2\pi a_2 \cos(\pi x) \\
 - \,\,\ds\pi \sum_{k\geq2} \left(\left(k+1\right)a_{k+1} - \left(k-1\right)a_{k-1}\right) \cos(k\pi x)\,. 
\end{multline}
\begin{multline}
\label{eq:fourrier3}
\cos(\varphi^c) \partial^3_x\varphi^w = - \pi^2 a_1 - 4 \pi^2 a_2 \cos(\pi x) \\
 - \,\,\ds\pi^2 \sum_{k\geq2} \left(\left(k+1\right)^2a_{k+1} + \left(k-1\right)^2a_{k-1}\right) \cos(k\pi x) 
\end{multline}
\begin{multline}
\label{eq:fourrier4}
\sin(\varphi^c) \partial^4_x\varphi^w =   \pi^3 a_1 + 8 \pi^3 a_2 \cos(\pi x)  \\
 + \,\,\ds \pi^3 \sum_{k\geq2} \left(\left(k+1\right)^3a_{k+1} - \left(k-1\right)^3a_{k-1}\right) \cos(k\pi x)\,. 
\end{multline}
In addition, the Fourier series of the functions $ \cos(\varphi^c) \mu^w $, 
$ \sin(\varphi^c) \partial_x\mu^w $, 
$ \cos(\varphi^c) \partial_x^2\mu^w $ and
$ \sin(\varphi^c) \partial_x^3\mu^w $ are
involved in the linearized equations \eqref{eq:linearized_2} and \eqref{eq:mu_linearize}. They
are obtained by replacing the function $\partial^k_x \varphi^w$ by $\partial_x^{k-1}\mu^w$ for $k=1,\cdots,4$ and the coefficients $a_i$ by $c_i$ for $i\geq0$, in formulas \eqref{eq:fourrier1}--\eqref{eq:fourrier4}.
We obtain that the Fourier coefficients for $k\geq 2$ are linked via the following relation
\begin{equation}
\label{eq:Meca_Fourrier1}
\dfrac d{dt}\left( M_1 A \right) = P L^c F  \left( M_2  A + \deg(F;\mu^c)   M_3 C\right)\,,
\end{equation}
where $A$ is the infinite column vector $
A=(a_k)_{ k\geq 2 }$,  $C = (c_k)_{ k\geq 2}$, 
and the   $M_i$, $i=1,2,3$, are upper triangular infinite matrices 
\[
M_1 =
\begin{pmatrix}
p_1(2) & 0 & p_2(2) & 0 & & &\\
 & \ddots &   \ddots & \ddots &  & & \\
& 0 &  p_1(k) &   0 &  p_2(k)  & 0 &  \\
 & &  & \ddots&  \ddots  & \ddots&  \\
\end{pmatrix}\,,
\]
$$ 
M_2 =
\begin{pmatrix}
q_1(2) & 0 & q_2(2) & 0 & & &\\
 & \ddots &   \ddots & \ddots &  & & \\
& 0  &  q_1(k) &   0 &  q_2(k)  & 0 &  \\
 & &  & \ddots&  \ddots  & \ddots&  \\
\end{pmatrix}\,,
$$
$$
M_3 =
\begin{pmatrix}
r(2) & 0 & -r(2) & 0 & & &\\
 & \ddots &   \ddots & \ddots &  & & \\
 & 0 &  r(k) &   0 &  -r(k)  & 0 &  \\
 &  &  & \ddots&  \ddots  & \ddots&  \\
\end{pmatrix}\,.
$$ 
The polynomials $p_i,q_i$, $i=1,2$ and $r$ are expressed for $k\geq2$ as
\begin{align*}
&p_1(k) = \pi \left(2+ k\right)\,,\quad p_2(k) = \pi  \left(2- k\right)\, ,  \\
&q_1(k) = -k^3 -k^2 + \left(3-\nu \right) k + 2(1-\nu)\, , \quad q_2(k) = k^3 -k^2  -\left(3-\nu\right) k + 2(1-\nu)\, , \\
&r(k) =\pi \left(1-\nu\right) k  \left(k^2 -4\right)\, .
\end{align*}
Since   $M_1$ does not depend on time, the set of equations \eqref{eq:Meca_Fourrier1} becomes
\begin{equation}
\label{eq:Meca_Fourrier2}
\dfrac d{dt} A   =  P L^c F  \left( M_1^{-1}M_2  A +  \deg(F;\mu^c)    M_1^{-1}M_3 C\right)\,.
\end{equation}

\begin{remark}[Competition between mechanics and expansion]
We notice already some interesting properties.
First, the roots of polynomial $q_1$ are $-2$, $\frac{1}{2} \pm \frac{1}{2} \sqrt{5-4\nu}$.
Thereby, for $k\geq2$, $q_1(k) \leq0$ since $\nu \in (0,1)$. 
Consequently, all the diagonal coefficients of  $M_1^{-1}M_2$, which are equal to $\frac{q_1(k)}{p_1(k)}$, are negative. It underlines the fact that the mechanics stabilizes the spherical solution.   
On the other hand, $M_3$ represents the coupling with heterogeneous growth. Since the roots of polynomial $r$ are $-2$, $0$ and $2$, then $\frac{r(k)}{p_1(k)} \geq 0$ for $k\geq2$, heterogeneous growth is clearly a destabilizing phenomena in this setting. The present analysis enables to quantify this competition.
\end{remark}

By using the linearized equation \eqref{eq:mu_linearize}, we obtain
\begin{equation}
\label{eq:fourier_mu}
C = N_1^{-1}N_2 A\,,
\end{equation}
where 
$$
N_1 = 
\begin{pmatrix}
s_1(2) & 0 & s_2(2) & 0  & & &\\
 & \ddots &   \ddots & \ddots &  & & \\
& 0 &  s_1(k) &   0 &  s_2(k)  & 0 &  \\
 &  &  & \ddots&  \ddots  & \ddots&  \\
\end{pmatrix}\, , 
$$ 
and
$$
N_2 =
\begin{pmatrix}
t_1(2) & 0 & t_2(2) & 0 & & &\\
 & \ddots &   \ddots & \ddots &  & & \\
& 0 &  t_1(k) &   0 &  t_2(k)  & 0 &  \\
 & &  & \ddots&  \ddots  & \ddots&  \\
\end{pmatrix}\,.
$$ 
The coefficients $s_i(k)$ and $t_i(k)$, $i=1,2$ are expressed for $k\geq2$  as
\begin{align*} 
& s_1(k) =\pi \left( k-\sigma^c{k}^{2}+\sigma^c{k}^{3}\right) \, , \quad s_2(k) = - \pi \left( k+  \sigma^c{k}^{2}+  \sigma^c {k}^{3}\right)\, , \\
&  t_1(k) =  k+1\, , \quad 
t_2(k) =- \left(k-1\right)\, .
\end{align*}
By replacing $C$ in \eqref{eq:Meca_Fourrier2} by its expression in \eqref{eq:fourier_mu}, we get eventually
\begin{equation}
\label{eq:Meca_Fourrier3}
\dfrac d{dt} A  =  \vecM \,A \,,
\end{equation}
where 
\begin{equation}
\label{eq:Meca_Fourrier2bis}
\vecM  = PL^c F \left(  M_1^{-1}M_2 + \deg(F;\mu^c)  M_1^{-1}M_3N_1^{-1}N_2\right)\,.
\end{equation}
All the matrices appearing in the above expression  are in fact upper triangular. The diagonal coefficients of $\vecM$ are denoted by $\lambda_k$. For $k\geq2$, we have
$$
\lambda_k =\frac{PL^c F }{p_1(k)}\left( q_1(k) + \deg(F;\mu^c)\frac{r(k)t_1(k)}{s_1(k)}\right)  =  P L^c  F \frac{N(t,k)}{D(t,k)}\, ,
$$
where the numerator is given by
$$
 N(t,k) = n_0(t) + n_1(t)  k + n_2(t)  k^2 + n_3(t)  k^3 +n_4(t)  k^4  \,, $$
with
\begin{align*}
n_0(t) &=  (1-\nu )\left[1-2 \deg(F;\mu^c(t))\right] \,,\\
n_1(t)  &= \,-(1-\nu)\left[\sigma^c(t) + \deg(F;\mu^c(t))\right] +1 \,,\\
n_2(t)  &=  -\nu\sigma^c(t) -1 + \left( 1-\nu \right) \deg(F;\mu^c(t)) \,,\\
n_3(t) &= 2\sigma^c(t)  \,,\\
n_4(t) &= -\sigma^c(t) \,.
\end{align*}
and the denominator   is given by
$$
D(t,k) = \pi\left(\sigma^c(t) k^2 - \sigma^c(t) k + 1\right)\geq 0\,.
$$
%
%

\section{Numerical results}
\label{Sec:Numerics}

The numerical scheme is briefly introduced, and numerical results are discussed in the 3D case. In particular, we investigate the linear stability of the spherical shape, and we explore the behaviour of the system in the nonlinear regime. 
\bigskip



We apply standard finite difference methods for transport-diffusion equations. 
We opt for a piecewise linear approximation of the function $\varphi$ on the regular grid $0 = x_0< x_1< \dots< x_{m+1} = 1$ with space step $x_{i+1} - x_i = \Delta x$. The functions $\partial_x \varphi$ and  $\mu$ are approximated by piecewise constant functions on the shifted grid $\Delta x/2 = y_0 < y_1 < \dots < y_m = 1 -  \Delta x/2$. The transport parts, respectively in \eqref{eq:2D_new_Phi}, and \eqref{eq:3D_new_Phi} are approximated by a first-order upwind scheme. As for the diffusion terms, we use a semi-implicit approximation. 
We have paid much attention to the discretization at the boundary in order to remove the possible singularities. For instance, the approximation of the integral $\int_0^x \cos \varphi(t,x')\, \mbox{d}x'$ \eqref{eq:Icos}, defined on the shifted grid, is computed using the piecewise linear approximation of $\varphi$.

We checked carefully that our numerical results are in agreement with the linear stability analysis.  
\bigskip

The initial condition $\varphi(0,x)$ is a perturbation of the spherical shape $\varphi(x) = \pi x$ having ten Fourier modes with random coefficients. 

The coupling function $F$ is homogeneous, $F(\mu) = \mu^d$, so that $\deg(F;\mu) \equiv d$. We shall discuss the influence of $d$, the degree of nonlinearity, and  $\sigma$ the membrane diffusion coefficient. An important point is that $\sigma$ depends on $t$ due to the increase in length (recall that $\sigma(t)=\frac{\gamma \pi^2}{\alpha L(t)^2}$). 
\paragraph{Quasi-stationary expansion.}

\begin{figure} 
\begin{center}
\includegraphics[width=0.45\linewidth]{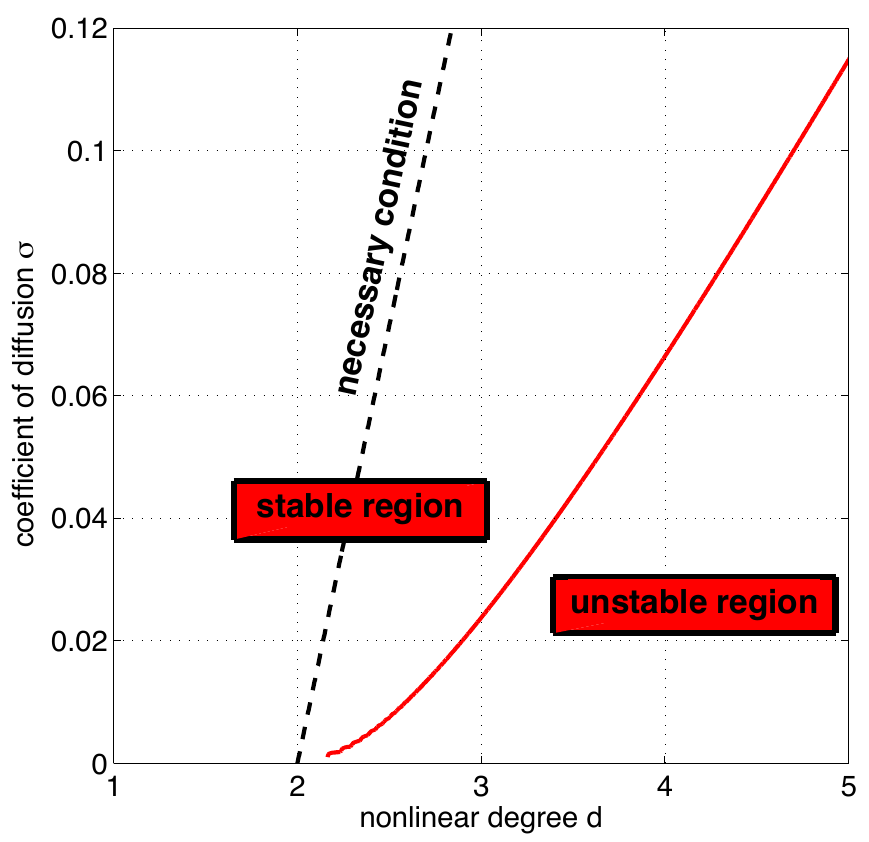}{\footnotesize(a)}\;
 \includegraphics[width=0.45\linewidth]{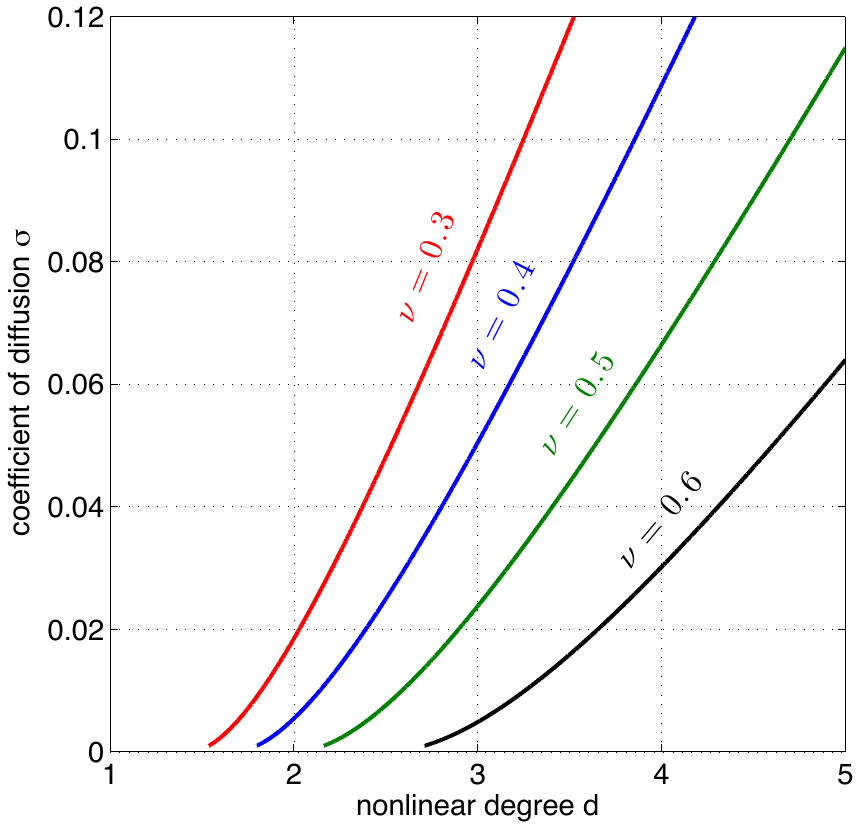}
{\footnotesize(b)}
\caption{\small \label{Fig:condition_L}(a) Stability region in the space $(d,\sigma)$ for $\nu = 0.5$. The red curve is the set of points such as the difference between the two greatest roots of polynomial \eqref{eigenval} are equal to $1$. 
For the sake of comparison, the necessary condition derived in \eqref{eq:necessary} is plotted in dashed line. (b) Influence of the Poisson's ratio $\nu$.}
\end{center}
\end{figure}

%

Firstly, we suppose that cell expansion is slow enough to consider that $L(t)$ is approximately constant. 

Unsurprisingly, the instability region in the $(d,\sigma)$ space lies below a curve which is increasing with respect to $d$ (Figure \ref{Fig:condition_L}(a)). Therefore, the instability is favoured as $d$ increases or $\sigma$ decreases. As a by-product, instability is enhanced as the cell is expanding, since $\sigma$ scales as the inverse of $L(t)^2$. On the other hand, for a given $(d,\sigma)$, the instability region gets larger when $\nu$ increases  (Figure \ref{Fig:condition_L}(b)). 

\begin{figure} 
\begin{center}
\includegraphics[width = 0.45\linewidth]{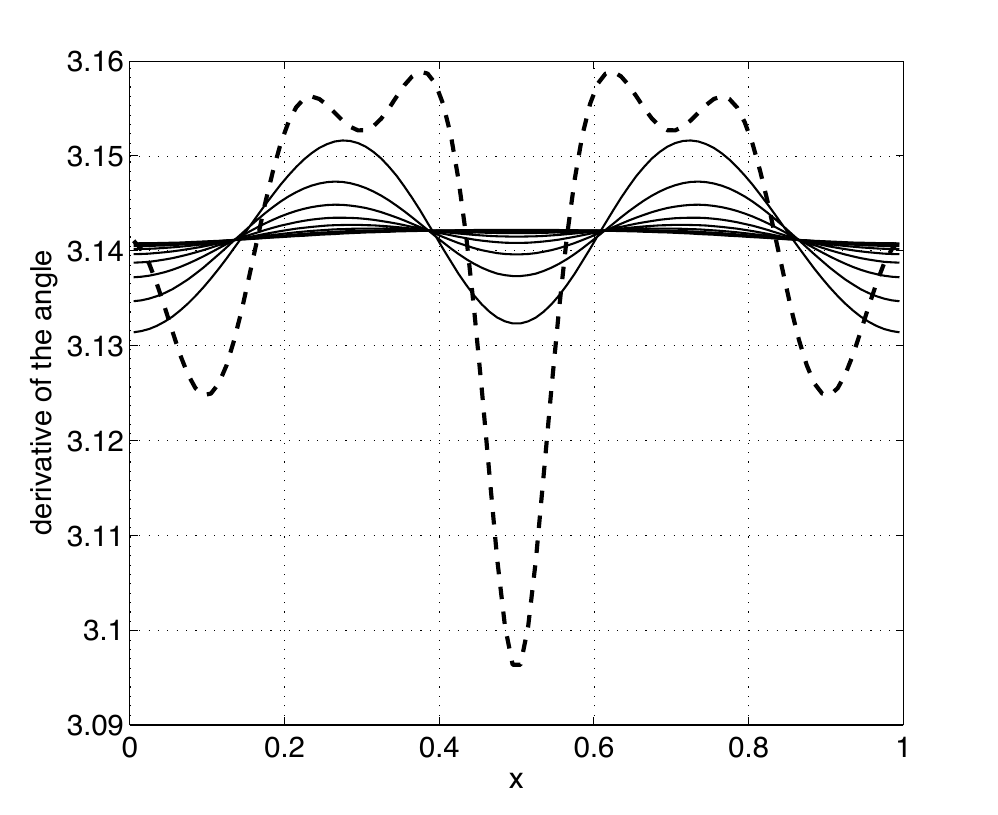}
\caption{\small\label{Fig:3D_stable_pattern}  Evolution of function $\partial_x \phi$ for a stable set of parameters $\sigma = 10^{-1}$ and $d=4$. Initial shape $ \partial_x \phi(0, x)$ is plotted in dashed line.}
\end{center}
\end{figure}

\begin{figure} 
\begin{center}
\includegraphics[width = 0.45\linewidth]{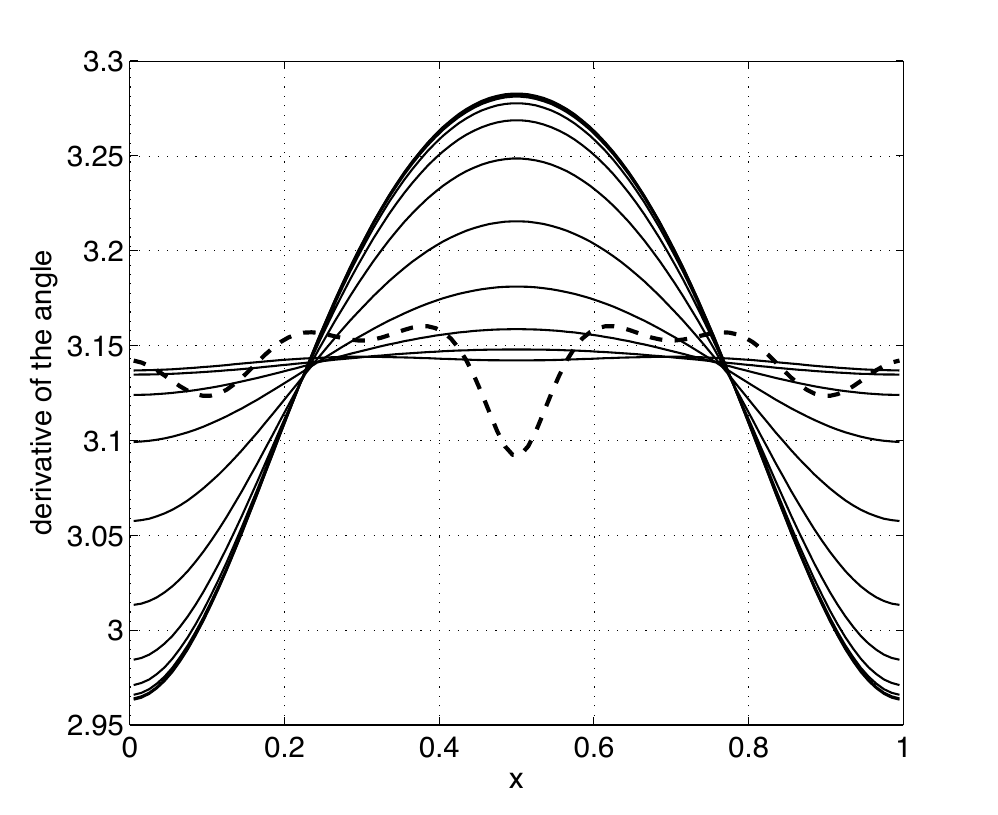}{\footnotesize(a)}\; \includegraphics[width = 0.45\linewidth]{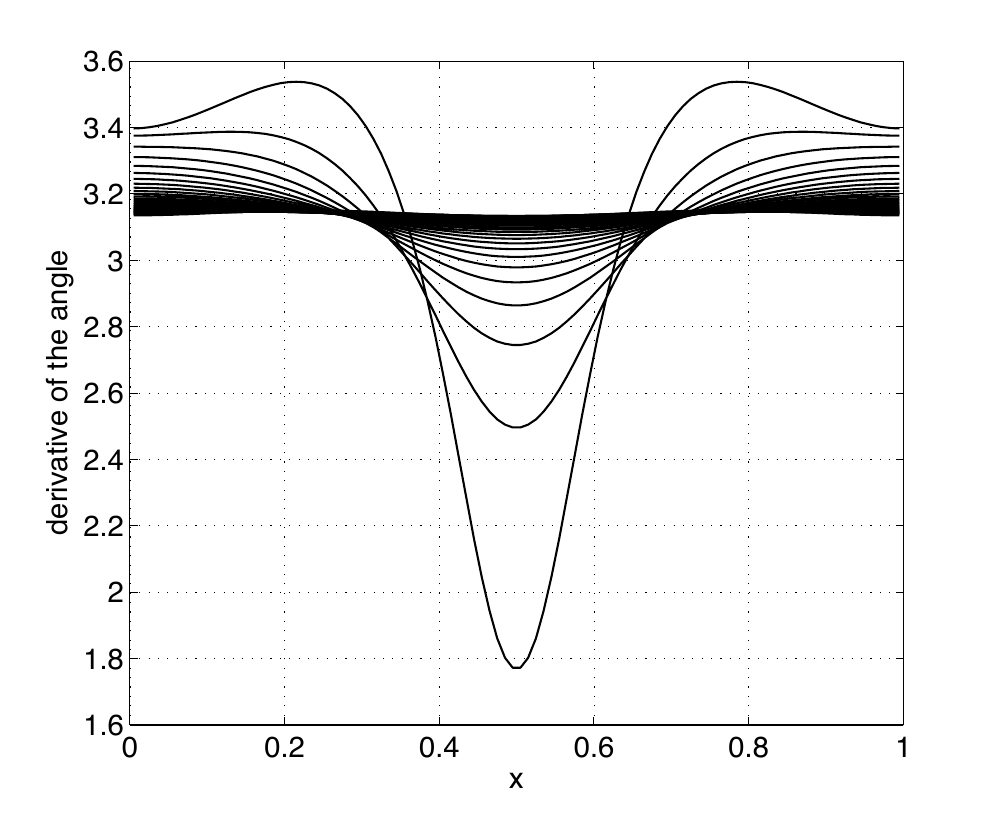}{\footnotesize(b)}
\caption{\small\label{Fig:3D_unstable_pattern} Evolution of  $\partial_x \phi$ for an unstable set of parameters $\sigma = 5.10^{-2}$ and $d=4$, and two different initial conditions. Initial shapes $\partial_x \phi(0, x)$ are plotted in dashed line.}
\end{center}
\end{figure}

We present in Figure \ref{Fig:3D_stable_pattern} a numerical simulation for a set of parameters for which the spherical shape is stable. We observe that the initial perturbation is quickly relaxed towards the constant state $\partial_x \varphi \equiv \pi$. We present in Figures \ref{Fig:3D_unstable_pattern} two numerical simulations for the same set of unstable parameters, but with different initial conditions. This is to illustrate various possible nonlinear behaviors.  We observed at least two possible profiles, both in the two-dimensional case and the three-dimensional case: i) $\partial_x\phi$ is above the mean for intermediate $x$ (Figure \ref{Fig:3D_unstable_pattern}(a)), or ii) $\partial_x\phi$ is above the mean for extremal $x$ (Figure \ref{Fig:3D_unstable_pattern}(b)). To explain the selection of two possible modes, we opt for a simplified description  in 3D. We consider an inflating ellipsoid, the  axes of which are denoted by  $a(t), b(t), c(t)$. By axisymmetry, we assume w.l.o.g. that $a(t) = b(t)$. We assume the following rule for the evolution of the curve: the rate of expansion of the axes is equal to the Gaussian curvature at the corresponding tips. Therefore we write
\begin{equation}\label{eq:ellipsoid}
\left\{\begin{array}{l}
\dfrac {d a(t)}{dt} = K = \dfrac{a(t)^2}{b(t)^2 c(t)^2} = \dfrac{1}{c(t)^2} \medskip \\ 
\dfrac {d c(t)}{dt} = K =  \dfrac{c(t)^2}{a(t)^2b(t)^2} = \dfrac{c(t)^2}{a(t)^4}\, .
\end{array}\right.
\end{equation} 
We notice that the quantity $a(t)^{-3} - c(t)^{-3} $ is constant for $t>0$. We denote $\delta = a(0)^{-3} - c(0)^{-3}$. The first equation in \eqref{eq:ellipsoid} becomes 
\[ \dfrac {d a(t)}{dt} =   \left( \dfrac{1}{a(t)^{3}} - \delta \right)^{2/3}\, . \]
The asymptotic behaviour depends on the sign of $\delta$. If $\delta>0$, {\it i.e.} $c$ is the major axis of the ellipsoid, then $a(t)$ converges asymptotically towards the globally stable equilibrium $a_\infty = \delta^{-1/3}$, whereas $c$ blows-up in finite time. On the contrary, if $\delta < 0$, {\it i.e.} $c$ is the minor axis of the ellipsoid then $a(t)$ grows linearly with time. On the other hand, the second equation in \eqref{eq:ellipsoid} becomes 
\[ \dfrac {d c(t)}{dt} =  c(t)^2 \left( \dfrac{1}{c(t)^{3}} + \delta \right)^{4/3}\, . \] 
therefore $c(t)$ converges asymptotically towards the globally stable equilibrium $c_\infty = (-\delta)^{-1/3}$. 
As a conclusion, depending on the relative length of the axes initially, the ellipsoid will evolves towards a cigar shape, or a flat expanding shape. Of course, all aspects of the mechanics have been neglected in this toy model.

\paragraph{Inflating system.}

\begin{figure} 
\begin{center}
\includegraphics[width = 0.45\linewidth]{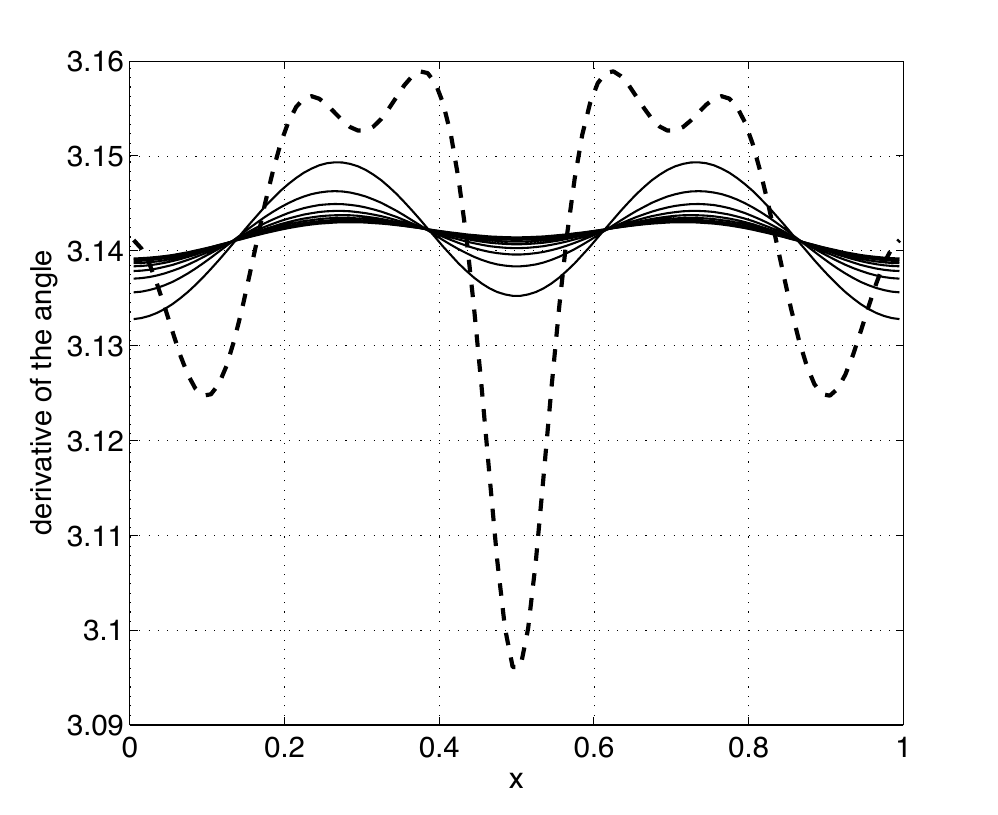}\; 
\includegraphics[width = 0.45\linewidth]{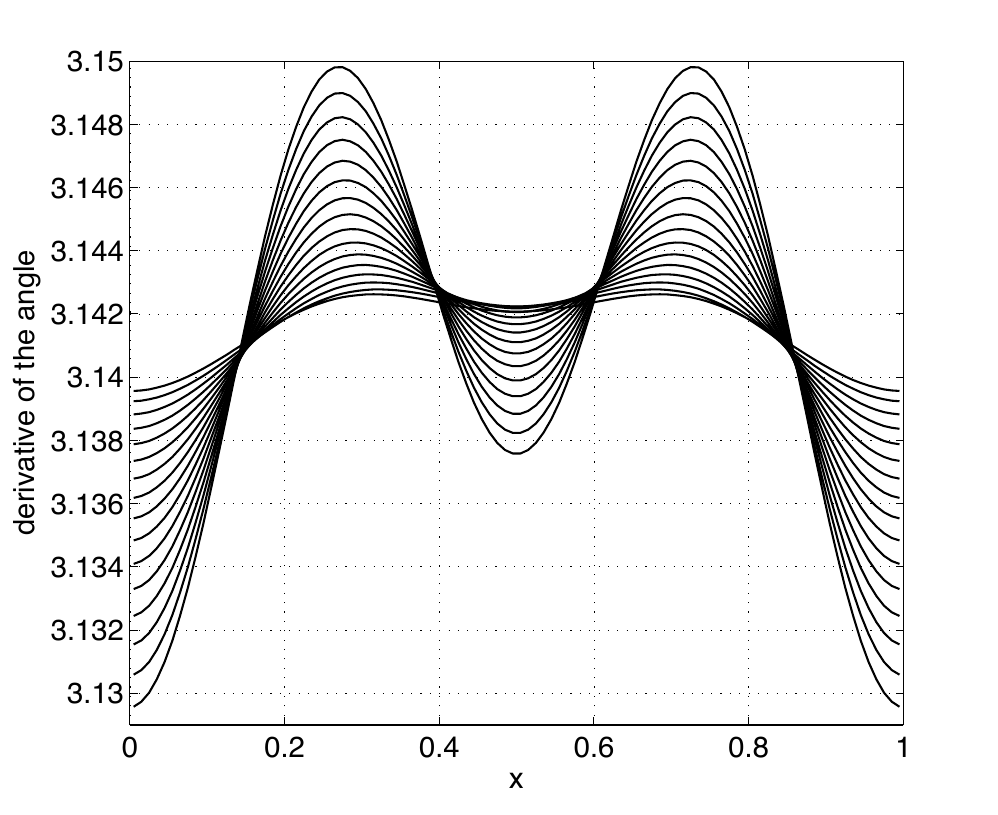}
\caption{\small\label{Fig:3D_stable_then_unstable01} Evolution of function $\partial_x \phi$ for the same set of parameters as in Figure \ref{Fig:3D_stable_pattern}, but for an increasing $L(t)$. The initial perturbation is damped (Left), but some instability arises as $L(t)$ gets larger after some time (Right).}
\end{center}
\end{figure}

Then, we investigate the evolution of the system in the case where $L(t)$ is free to increase (when the rate of expansion is comparable to the relaxation time of the system). We expect that the system becomes more likely unstable as $L(t)$ increases (see Figure \ref{Fig:condition_L}). We also expect more sophisticated patterns. In Figure \ref{Fig:3D_stable_then_unstable01}, we use the same set of parameters and initial condition as in \ref{Fig:3D_stable_pattern}. We observe the same dynamics for earlier time, but the spherical shape becomes progressively unstable, as expected.

\begin{figure}
\begin{center}
\hspace{-1cm}
\includegraphics[width = 0.55\linewidth]{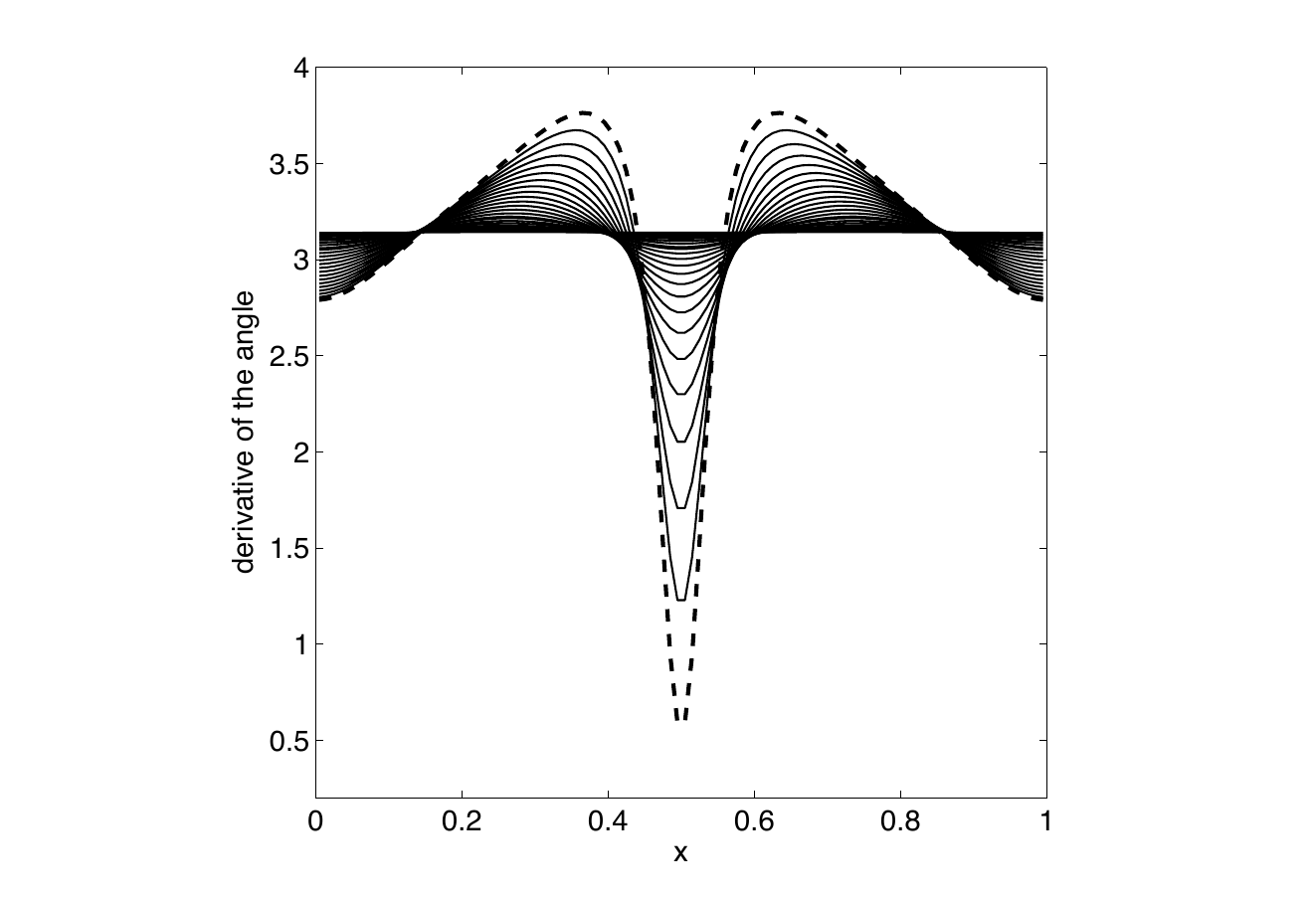}
\includegraphics[width = 0.37\linewidth]{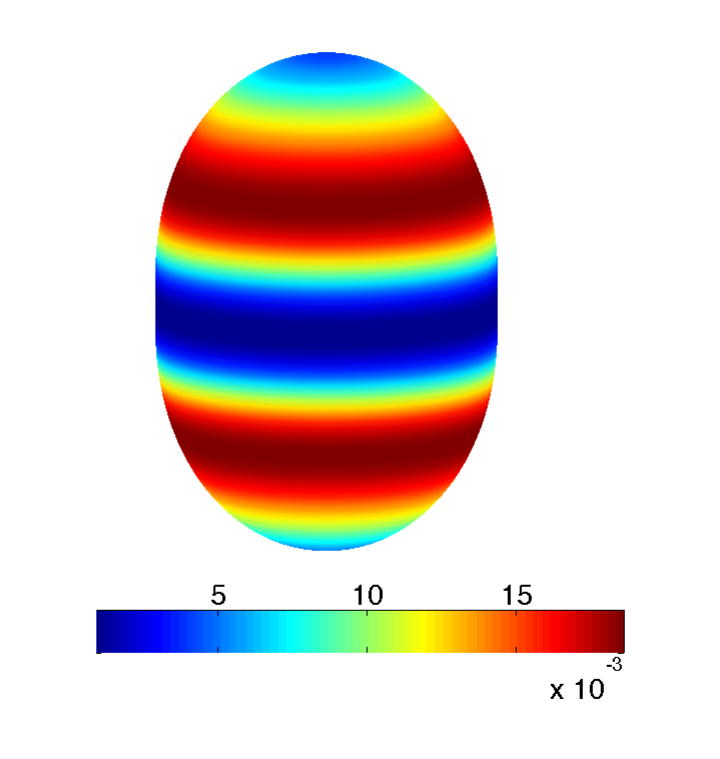}
%
\hspace{-0.8cm}

\includegraphics[width = 0.55\linewidth]{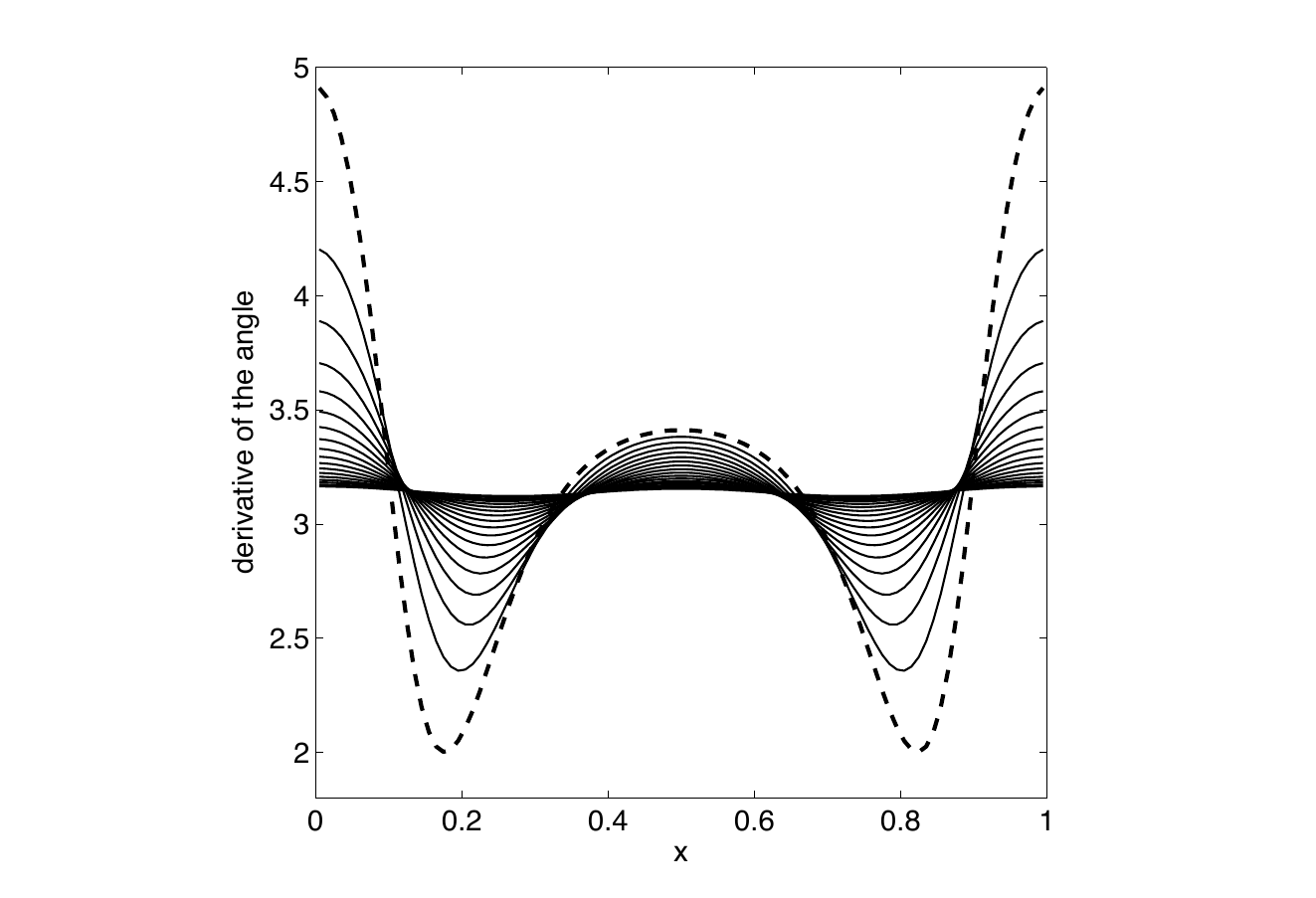}
\includegraphics[width = 0.37\linewidth]{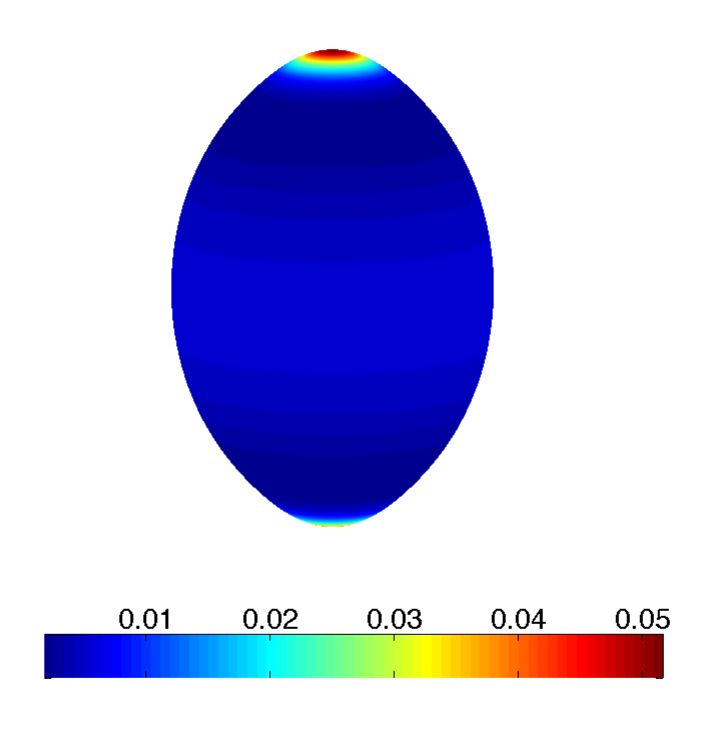}
\caption{\small\label{Fig:unstable01}Evolution of function $\partial_x \phi$ for an unstable set of parameters $\sigma = 5.10^{-2}$ and $d=4$, and two different initial conditions. Final shape plotted in dash line on the left corresponds to the surface plotted on the right. The color (online) represents the quantity of the density of materials available on the membrane.}
\end{center}
\end{figure}

Finally, we present in Figure \ref{Fig:unstable01} two possible behaviours in the nonlinear, inflating regime, as in Figure \ref{Fig:3D_unstable_pattern}. We also show the reconstruction of the three-dimensional axisymmetric cell as a result of the numerical simulations. We observe clearly the cigar shape up in Figure \ref{Fig:unstable01}, and the rod shape below in Figure \ref{Fig:unstable01}. The latter is typical of fission yeast growth. 


\vspace{1cm}

\noindent\textbf{Acknowledgements:}
Laetitia Giraldi was funded by the french ANR project MODPOL ANR-11-JS01-003-01.
This work was initiated within the framework of the LABEX MILYON (ANR-10-LABX-0070) of Université de Lyon, within the program "Investissements d'Avenir" (ANR-11-IDEX-0007) operated by the French National Research Agency (ANR).
This work was partially supported by the Labex LMH through the grant  ANR-11-LABX-0056-LMH in the "Investissements d'Avenir''.
The authors express their gratitude to P. Vigneaux for invaluable discussions.
We also thank A. Boudaoud, M. Piel and J. C. Sikorav and R. Voituriez for fruitful discussions at the early stage of this project. Finally, we thanks M. Bacha for its interest in this subject.

\bibliographystyle{plain}
\bibliography{Croissance_Cell}

\newpage

\appendix
\section{Kinematics}
\label{sec:kinematics}
In this appendix, we give more details on the way to derive equation \eqref{eq:r_final} governing the evolution of the generatrix curve $\mathcal{C}_t$ pushed by a vector fields, $\vecv(t,\cdot): \mathcal{C}_t \to \mathbb{R}^2$.\\

Let $t$ be a positive real and $s \in [0,L(t)]$, we express $\frac{\partial }{\partial t}\left(\frac{\partial r(t,s(t))}{\partial s}\right)$ as 

\begin{align}
\lefteqn{\frac{\partial }{\partial t}\left(\frac{\partial r(t,s(t))}{\partial s}\right) =}\nonumber \\
&\lim_{\mbox{dt},\mbox{d}\tilde{s}\, \to 0}\frac{1}{\mbox{dt}}  \left( \frac{r(t+\mbox{dt},\tilde{s}+\mbox{d}\tilde{s}) -r(t+\mbox{dt},\tilde{s})}{\mbox{d}\tilde{s}} - \frac{\partial r(t,s)}{\partial s}\right)\,,\label{eq_taux_r}
\end{align}

where $\mbox{dt}$ (resp. $\mbox{d}\tilde{s}$ and then $\mbox{ds}$) is an infinitesimal time (resp. length). We denote by $s$ the curvilinear abscissa at time $t$ and by $\tilde{s}$ the one at time $t+\mbox{dt}$.

However, we have (see Figure \ref{Fig:eq_evolution})

\begin{align}
\label{eq:r_infinitesimal}
r(t+\mbox{dt},\tilde{s}) &= r(t,s) + \left(\vecv\cdot \vece_1 \right)\mbox{dt}+o(\mbox{dt}) \,,\nonumber\\
r(t+\mbox{dt},\tilde{s}+\mbox{d}\tilde{s}) & = r(t,s+\mbox{ds}) + \left(\left(\vecv+\mbox{dv}\right)\cdot\vece_1\right) \mbox{dt}+o(\mbox{dt}) \,,
\end{align}
here the vectors $\mbox{dt} \vecv$  (resp. $\mbox{dt}\left(\vecv+\mbox{dv}\right)$) represents  the vector fields $\vecv(t,\cdot)$ at $\left(r(t,s),z(t,s)\right))$ (resp. at $\left(r(t,s+\mbox{ds}),z(t,s+\mbox{ds})\right))$ as depicted in Figure \ref{Fig:eq_evolution}. 
\\

Using the relations \eqref{eq:r_infinitesimal}, we get
\begin{equation}
\label{eq_term}
\frac{r(t+\mbox{dt},\tilde{s}+\mbox{d}\tilde{s}) -r(t+\mbox{dt},\tilde{s})}{\mbox{d}\tilde{s}} = \left(\frac{\partial r(t,s)}{\partial s} + \frac{\left(\mbox{dv} \cdot \vece_1\right)\mbox{dt}}{\mbox{ds}} \right)\frac{\mbox{ds}}{\mbox{d}\tilde{s}}+o(\mbox{dt}) \,.\\
\end{equation}


Notice that the derivative of the vector field $\vecv$ with respect to $s$ is 
\begin{equation}
\label{eq:partial_s_v}
 \frac{\partial \vecv(t,s)}{\partial s} \cdot \vece_1 =\frac{\left(\mbox{dv} \cdot \vece_1\right)}{\mbox{ds}}\,.
\end{equation}

Moreover, the infinitesimal quantity $\mbox{d}\tilde{s}$ reads (see Fig )
\begin{align*}
\label{eq_ds}
\mbox{d}\tilde{s} &=\mbox{ds} - \mbox{dt} \left(\vecv\cdot \vectau\right) + \mbox{dt} \left(\left(\vecv+\mbox{dv} \right) \cdot \vectau\right)+  o(\mbox{dt})\,, \nonumber \\
&= \mbox{ds} + \mbox{dt}\mbox{dv} \cdot \vectau+  o(\mbox{dt})
\end{align*}

By dividing the previous equality by $\mbox{ds}$ and inverting it, we get
\begin{equation}
\label{eq:ds_sur_ds}
\frac{\mbox{ds}}{\mbox{d}\tilde{s}} = 1 - \mbox{dt} \frac{\partial \vecv(t,s)}{\partial s} \cdot \vectau +  o(\mbox{dt})\,.
\end{equation}

By plugging \eqref{eq:partial_s_v} and \eqref{eq:ds_sur_ds} into \eqref{eq_term}, we obtain
\begin{align}
\label{eq_term_2}
\frac{r(t+\mbox{dt},\tilde{s}+\mbox{d}\tilde{s}) -r(t+\mbox{dt},\tilde{s})}{\mbox{d}\tilde{s}} =& \frac{\partial r(t,s)}{\partial s} +  \mbox{dt} \frac{\partial \vecv(t,s)}{\partial s} \cdot \vece_1 \nonumber \\
& -\mbox{dt} \frac{\partial r(t,s)}{\partial s} \frac{\partial \vecv(t,s)}{\partial s} \cdot \vectau + o(\mbox{dt})\,.
\end{align}

By plugging \eqref{eq_term_2} into \eqref{eq_taux_r} and by expanding $\frac{\partial }{\partial t}\left(\frac{\partial r(t,s(t))}{\partial s}\right)$, we get the equation, governing the evolution of the generatrix curve for all $t>0$ and for all $s \in [0,L(t)]$,

\begin{equation*}
\label{eq:r_expand}
\frac{\partial }{\partial t}\left(\frac{\partial r(t,s)}{\partial s}\right)   + \frac{ds}{dt} \frac{\partial^2 r(t,s)}{\partial^2 s} = - \frac{\partial r(t,s)}{\partial s} \left(\frac{\partial \vecv(t,s)}{\partial s} \cdot \vectau(t,s)\right) + \frac{\partial \vecv(t,s)}{\partial s} \cdot \vece_1\,.
\end{equation*}

Finally, noticing $\ds\frac{ds}{dt} =\ds \int_0^s( \frac{\partial  \vecv(t,s')}{\partial s}) \cdot\vectau(t,s')\, \mbox{ds}'$, we obtain for all $t>0$ and for all $s \in [0,L(t)]$ the equation \eqref{eq:r_final}.


\section{Change of variables}
\label{Appendix:Lemma}


This section is devoted to some preliminary simplifications of the system of equations \eqref{eq:r_final}-\eqref{eq:velocity_stress}-\eqref{eq:Laplacian_Beltrami}. 
Firstly, we change the parametrization of the generatrix curve and derive its associated evolution equation. 
This leads to simplify the evolution equation \eqref{eq:r_final}.
Furthermore, through changing variable, we bring back the study of the system \eqref{eq:r_final}-\eqref{eq:velocity_stress}-\eqref{eq:Laplacian_Beltrami} defined on a varying-time domain into a fixed one. 
Finally, we combine the velocity vector fields equations \eqref{eq:velocity_stress} with the evolution one to get a condensed formula for the system of equations \eqref{eq:r_final}-\eqref{eq:velocity_stress}-\eqref{eq:Laplacian_Beltrami}.
\\


\noindent\textbf{Change of parametrization}\\
In what follows, instead to reconstruct the generatrix curve by using the system of coordinates $\left\{(s,r(t,s))\right\}$, we will use the parametrization given by $\left\{(s,\varphi(t,s))\right\}$. \\

By using relation \eqref{Change_Variable}, equation \eqref{eq:r_final} can be expressed as a partial differential equation on $(t,s)\mapsto \varphi(t,s)$.
Indeed, the left hand side of equation \eqref{eq:r_final} reads,
\begin{equation}
\label{eq:LHS_eq_phi}
\forall t >0\,, \forall s\in [0,L(t)]\,, \partial_t \left(\partial s r(t,s)\right) = - \sin(\varphi(t,s))\, \partial_t \varphi(t,s)\,.
\end{equation}

By using the expression of the two vector fields $\vecn$ and $\vectau$ (given in Fig \ref{Fig:surface_revolution}), we get 
\begin{equation}
\label{eq:v_change_var}
\vecv \cdot \vece_1 = v_n \sin(\varphi) + v_{\tau} \cos(\varphi)\,,\quad
\vecv \cdot \vectau = v_n \,\partial s \varphi + \partial s v_{\tau}\,.
\end{equation}

Consequently, by plugging the relations \eqref{eq:v_change_var} into the right hand side of equation \eqref{eq:r_final}  we obtain $\forall t >0\,, \forall s\in [0,L(t)]\,,$

\begin{align}
\lefteqn{\partial s \left( -\left(\int_0^s(\partial s\vecv) \cdot\vectau \right) \partial s r + \vecv\cdot\vece_1\right)=} \nonumber\\
&\sin(\varphi) \bigg(\left(\ds-\int_0^s v_n \partial s\varphi\right)\partial s \varphi - \partial s v_n\bigg) \label{eq:RHS_eq_phi}
\end{align}

Then, by dividing \eqref{eq:LHS_eq_phi} and \eqref{eq:RHS_eq_phi} by $\sin(\varphi)$, we have for all $t>0$ and $s\in [0,L(t)]$
 \begin{equation}
 \label{eq:phi_s}
 \partial_t \phi(t,s) = \left(\ds-\int_0^s v_n(t,s')  \partial s\varphi(t,s')\mbox{ds}' \right)\partial s \varphi(t,s) - \partial s v_n(t,s) \,,
 \end{equation}
where $v_n$ satisfies the equation \eqref{eq:vn}.
\\

Moreover, it derives from relation \eqref{Change_Variable} that for all $t>0$ and for all $s\in [0,L(t)]$,
 
\begin{align}
r(t,s) &= \int_0^s \cos(\varphi(t,s'))\mbox{ds}'\,,\nonumber\\
\kappa_s(t,s) &= \partial s\varphi\,, \nonumber \\
\kappa_{\theta}(t,s) &=\ds \frac{\sin(\varphi)}{\ds\int_0^s\cos(\varphi(t,s'))\mbox{ds}'}\,. \label{mu_change_para}
\end{align}

By using relations \eqref{mu_change_para}, equation \eqref{eq:mu} becomes 

\begin{align}
\lefteqn{
\beta \left(\frac{\partial s\varphi(t,s)\sin(\varphi(t,s))}{\ds\int_0^s\cos(\varphi(t,s'))\mbox{ds}'}\right)=}\nonumber\\ 
&-\gamma\frac{\ds\partial s\left(\int_0^s \cos(\varphi(t,s'))\mbox{ds}' \partial s \mu(t,s) \right) }{\ds\int_0^s \cos(\varphi(t,s'))\mbox{ds}'}+ \alpha\mu(t,s) 
\label{eq:mu_phi}
\end{align}

\noindent\textbf{Change of variables}\\
The system governing the dynamics of the cell membrane is posed on a time-varying domain $[0,L(t)]$. 
By setting $x:=\ds\frac{s}{L(t)}$, the system \eqref{eq:phi_s}-\eqref{eq:mu_phi} can be expressed on the fixed domain $[0,1]$.   
\\

In what follows, for all the function $f$ defined on $\mathbb{R}^+\times[0,L(t)]$,
we denote by $\tilde{f}$ the function defines on $\mathbb{R}^+\times[0,1]$ such that 
$$
\forall t>0\,, \quad \forall x\in [0,1]\,, \quad \tilde{f}(t,x) = f(t,xL(t))\,.
$$


By replacing $s$ by $xL(t)$ and using Definition \ref{eq:relation_tilde},  the equation \eqref{eq:phi_s} becomes for all $t>0$ and $\forall x\in[0,1]$,
\begin{align}
 \partial_t \tilde{\varphi}(t,x) =&\ds\Bigg( \frac{\partial_t L(t)}{L(t)} x - \left(\ds\int_0^x \tilde{v}_n(t,x')  \partial_x\tilde{\varphi}(t,x')\mbox{dx}' \right)\Bigg)\frac{\partial_x \tilde{\varphi}(t,x)}{L(t)} \nonumber\\
 &\ds\quad -\frac{ \partial_x \tilde{v}_n(t,x)}{L(t)} \label{eq:phi_x}\,,
 \end{align}
As well, the normal velocity field $\tilde{v}_n$ is equal to
\begin{equation}
 \label{eq:tilde_v_n}
 \tilde{v}_n(t,x) = \cot(\tilde{\varphi}(t,x)) \tilde{v}_t(t,x) + H(t,x)\,,
\end{equation}

where the function $H$ is defined on $\mathbb{R}^+\times [0,1]$ as
$$
H(t,x) = \frac{\tilde{\Psi}(t,x)}{\tilde{\kappa}_\theta(t,x)}\left(\tilde{\sigma}_\theta(t,x)-\nu \tilde{\sigma}_{s}(t,x)\right)\,.
$$

Moreover, since $v_{\tau}$ is solution to the ODE \eqref{eq:vt}, the tangential velocity field $\tilde{v}_t$ is the solution of the following ODE
 
\begin{equation}
 \partial_x\tilde{v}_t(t,x) - \Big(\partial_x\left( \tilde{\varphi}(t,x)\right)\, \cot(\tilde{\varphi}(t,x))\,\Big) \,\tilde{v}_t(t,x) = d(t,x) 
\label{eq:tilde_v_t}
\end{equation}
where the function $d$ is defined on $\mathbb{R}^+\times [0,1]$ as
$$
d(t,x) =L(t) \tilde{\Psi}(t,x) \bigg(\tilde{\sigma}_s(t,x) - \nu \tilde{\sigma}_\theta(t,x)\bigg)\, - L(t) \tilde{\kappa}_s(t,x) H(t,x)\,.
$$

Let us expand, by using the equalities \eqref{eq:tilde_v_n} and \eqref{eq:tilde_v_t},

\begin{align}
\lefteqn{\ds \int_0^x \tilde{v}_n(t,x') \partial_x\tilde{\varphi}(t,x')\mbox{dx}'=}\nonumber\\
 & -\tilde{v}_t 
+ \int_0^x \bigg(d(t,x') + H(t,x')\,\partial_x(\tilde{\varphi}(t,x'))\bigg) \mbox{dx}' \label{eq:rel_1}
\end{align}

However, by differentiating $\tilde{v}_n$ with respect to $x$, we get
\begin{equation}
\label{eq:rel_2}
\partial_x \tilde{v}_n = \partial_x \tilde{\varphi} \,\,\tilde{v}_t - \cot(\tilde{\varphi})\,\, d + \partial_x H
\end{equation}

Then, plugging \eqref{eq:rel_1}-\eqref{eq:rel_2} into \eqref{eq:phi_x}, the function $\tilde{\varphi}$ is solution to the following partial differential equation
\begin{align}
\partial_t \tilde{\varphi}(t,x) =& \ds\left(\frac{\partial_t L(t)}{L(t)} x - \int_0^x d(t,x) - \int_0^x H(t,x) \,\partial_x\tilde{\varphi}(t,x)\right)\partial_x \tilde{\varphi}(t,x)\nonumber \\
&\quad+ \cot(\tilde{\varphi}(t,x)) \,d(t,s) -\partial_x H(t,x)\,.\label{eq:evolution_Phi_3D}
\end{align}


Let us call $\hat{\mu} = \frac{\alpha L^2}{\beta \pi^2}\tilde{\mu}$.
Finally,  equation \eqref{eq:mu_phi} could be expressed  for all times $t$ on a fixed domain $[0,1]$ as

\begin{align}
\label{eq:3D_new_mu_App}
  - \dfrac{\sigma(t)}{\pi^2} \left( \partial^2_x \hat{\mu}  + \frac{\cos(\tilde{\varphi})}{\mbox{Icos}}\partial_x\hat{\mu} \right) + \hat{\mu}  = \frac1{\pi^2}  \partial_x\tilde{\varphi}\,   \frac{\sin(\tilde{\varphi})}{\mbox{Icos}}  \,,
\end{align}
with the boundary conditions $\mbox{Icos}(t,0)\partial_x \hat{\mu}(t,0)=\mbox{Icos}(t,1)\partial_x \hat{\mu}(t,1) = 0$ and where $\sigma(t)=\frac{\gamma \pi^2}{\alpha L(t)^2}$. 

\end{document}